

\magnification=1200
\documentstyle{amsppt}
\define\G{\Gamma}
\define\Z{{\Bbb Z}}
\define\Q{{\Bbb Q}}
\define\C{{\Bbb C}}
\define\R{{\Bbb R}}
\define\la{\langle}
\define\ra{\rangle}
\define\De{\Delta}
\define\T{\frak{T}}

\define\Mod {\frak{M}}
\define\CC {\Cal{C}}
\define\SS {\Cal{S}}
\define\ZZ {\Cal{Z}}
\define\ZZZ {\Cal{Z}_{\bullet}}
\define\GG {\Cal{G}}
\define\cA{\Cal{A}}

\define\Aut{\operatorname{Aut}}
\define\reg{\text{reg}}
\define\sing{\text{sing}}

\define\qc{\text{q.c.}}
\define\Bl{\operatorname{Bl}}

\define\conf{\operatorname{Conf}}
\define\aconf{\operatorname{Aconf}}
\define\supp{\operatorname{supp}}

\define\st{\text{st}}
\define\sst{\text{sst}}
\define\pr{\text{pr}}

\define\refer#1{(#1)}

\define\bowep{1}

\define\harera{3}
\define\harerb{4}
\define\harvey{6}
\define\konts{6}
\define\milpen{7}

\define\strebel{9}
\define\wolpert{10}

\newcount\headnumber
\newcount\labelnumber
\newcount\fignumber
\define\section{\global\advance\headnumber
by1\global\labelnumber=0{{\the\headnumber}.\ }}
\define\label{(\global\advance\labelnumber by1 \the\headnumber
.\the\labelnumber )\enspace}
\define\fig{\global\advance\fignumber by1{Fig.\ \the\fignumber} }

\NoBlackBoxes
\topmatter
\title
Cellular decompositions of compactified moduli spaces of pointed curves
\endtitle
\rightheadtext{Compactified moduli spaces}
\author
Eduard Looijenga
\endauthor
\address
Faculteit Wiskunde en Informatica,
Rijksuniversiteit Utrecht,
PO Box 80.010, 3508 TA Utrecht,
The Netherlands
\endaddress
\email
looijenga\@math.ruu.nl
\endemail

\abstract
To a closed connected oriented surface $S$ of genus $g$ and a nonempty finite
subset $P$ of $S$ is associated a simplicial complex (the arc complex) that
plays a basic r\^ ole in understanding the mapping class group of the pair
$(S,P)$. It is known that this arc complex contains in a natural way the
product
of the Teichm\"uller space of $(S,P)$ with an open simplex. In this paper we
give an interpretation for the whole arc complex and prove that it is a
quotient of a Deligne--Mumford extension of this Teichm\"uller space and a
closed simplex. We also describe a modification of the arc complex in the
spirit of Deligne--Mumford.
\endabstract
\endtopmatter

\head
Introduction
\endhead

Given a closed connected oriented differentiable surface $S$
of genus $g$ and a finite nonempty subset $P$ of $S$, then the mapping class
group $\G (S,P)$ of this pair is the group of isotopy classes of sense
preserving diffeomorphisms of $S$ that fix $S$ pointwise. Harer proved in a
series of papers some remarkable properties of the cohomology of the $\G (S,P)$
(see \cite{\harerb} for an overview). In this work a central r\^ole is played
by
various simplicial complexes with an action of an appropriate mapping class
group that have in common the property that stabilizers of simplices look like
simpler mapping class groups. The complex depends on the context, but in all
cases it can for a suitable pair $(S,P)$ be identified with a subcomplex of the
{\it arc complex} $A(S,P)$. That complex is defined as follows: the vertices of
$A(S,P)$ are ambient isotopy classes relative $P$ of embedded unoriented
nontrivial loops and arcs in $S$ that connect two (possibly identical) points
of $P$ and avoid all other points of $P$ (where a loop is considered trivial if
it bounds an open disk in $S-P$) and finitely many such vertices span a simplex
if we can respresent them by loops and arcs which do not meet outside $P$. We
note that there is a piecewise linear map $\lambda$ from $A(S,P)$ to the
simplex
$\De _P$ spanned by $P$ characterized by the property that it sends a vertex
represented by an arc (resp.\ a loop) to the barycenter of the $1$-simplex of
$\Delta _P$ spanned by its end points (resp.\ the vertex of $\Delta _P$
representing the base point).

An important property of this complex  is that its interior can be identified
with the product of the Teichm\"uller space $\T (S,P)$ of the pair $(S,P)$
(i.e., the space of isotopy classes relative $P$ of conformal structures on
$S$)
and the open simplex $\Delta _P^{\circ}$. We may therefore regard $A(S,P)$ as
an
extension of $\T (S,P)\times \Delta _P ^{\circ}$. In the applications alluded
to
there was no apparent need to know what this extension actually represents, and
that may have been the reason that question received little attention. (An
exception is the paper by Bowditch and Epstein \cite{\bowep} about which we
shall say more below.)
The situation changed with Kontsevich's work on a conjecture of Witten\
\cite{\konts}, where it became
essential to interpret the part of $A(S,P)$ of lying over $\Delta _P^{\circ}$.
In this article  Kontsevich states the answer but omits a proof.  The present
paper grew out the desire to supply one and one of our main results now
interprets all of $A(S,P)$ in terms of the Deligne--Mumford compactification of
the moduli space $\Mod _g^P:=\G (S,P)\backslash\T (S,P)$. For a precise
statement we refer to theorem \refer{8.6}. Suffice here to say that for every
nonempty
subset $Q$ of $P$ we describe a quotient space $K_Q\Mod _g^P$ of the
Deligne--Mumford compactification of $\Mod _g^P$ and for every inclusion
$Q\subset Q'$ a quotient map $K_Q\Mod _g^P\to K_{Q'}\Mod _g^P$ such that the
geometric realization of the associated simplicial space over $\Delta _P$ can
be
identified with the orbit space $\G (S,P)\backslash A(S,P)$. In particular, $\G
(S,P)\backslash A(S,P)$ is a quotient of the product of the  Deligne--Mumford
compactification and $\Delta _P$. We suspect that the compactifications
$K_Q\Mod _g^P$ and the maps between them can be constructed in the
category of projective varieties and morphisms so that $\G (S,P)\backslash
A(S,P)$ becomes the geometric realization of a simplicial object in this
category. We state the relevant conjectures in \refer{3.3}

An intermediate result of our proof is a combinatorial description \refer{11.5}
of (a
thickened version of) the Deligne--Mumford compactification. More precisely, we
equivariantly blow up $A(S,P)$ in a certain manner over its boundary (in the
PL-category) to get a cell complex of which the orbit space naturally maps to
$\overline{\Mod}{}_g^P\times\Delta _P$ with fibers products of simplices (or
finite quotients thereof). This description may be helpful in determining which
of the cohomology classes that Kontsevich introduced in $\Mod _g^P$ extend to
$\overline{\Mod}{}_g^P$. A paper by Milgram--Penner \cite{\milpen} alludes to a
combinatorial construction of the Deligne--Mumford compactification (for the
case that $P$ is a singleton), but it is not clear to us whether what these
authors have in mind coincides with our construction.

The article by Epstein and Bowditch mentioned above came to our attention after
this paper was essentially completed. It also gives an interpretation of the
arc
complex, but in this it differs from ours in two respects. First, it takes the
hyperbolic point of view (which gives rise to a different embedding of
thickened
Teichm\"uller space in the arc complex) and second, our description is solely
in
terms of the Deligne--Mumford compactification. (For these reasons it is not
clear to us whether it could take care of Kontsevich's assertion.) We adopted
their term {\it arc complex} and we adapted our notation a little in order to
avoid too blatant clashes with theirs.

\smallskip
The plan of the paper is as follows. The first seven sections are intended to
have to some extent the characteristics of a review paper and were written with
a nonexpert
reader in mind. Yet they do contain results that we have not found in the
literature. In the first section we collect facts about the Teichm\"uller
spaces. The next two sections deal with certain extensions of them: we describe
a boundary for Teichm\"uller spaces in the spirit of Harvey based on the
Deligne--Mumford compactification and we introduce the quotients of the
Deligne--Mumford compactification alluded to above. In section 4 we discuss
some properties of the complex $A(S,P)$. The next two
sections we introduce metrized ribbon graphs and explicate the relationship
between this notion and the complex $A(S,P)$. In section 7 we invoke the
fundamental results of Strebel, culminating in theorem \refer{7.5}. The
subsequent sections are of more technical nature. In section 8 we describe the
geometric objects that are parametrized by the points of $A(S,P)$. Our first
main theorem \refer{8.6} is also stated there, but its proof is postponed to
the last section. The remainder of the paper is mostly concerned with the
combinatorial
versions of notions related to the Deligne--Mumford compactification. In
section
9 we introduce stable ribbon graphs of which we claim that it is the
combinatorial analog of the notion of a stable curve. This is justified in
section 10, where we show that a metrized stable ribbon graph can be obtained
as
the limit of a one-parameter family of ordinary metrized ribbon graphs. In the
final section 11 we construct the modification $A(S,P)$ mentioned above and
prove our second main theorem \refer{11.5}, namely that this modification is
essentially a
thickened Deligne--Mumford extension of $\T (S,P)$. Once this has been
established, the proof of our first main theorem is easily completed.

\remark{Acknowledgements}
I thank K.\ Strebel for help with the proof of \refer{7.5},  S.\ Wolpert and
J.\ Koll\'ar for correspondence regarding \refer{3.3} and A.J.\ de Jong for the
observation mentioned in \refer{3.2}.
\endremark

\smallskip
Throughout this paper $S$ stands for a compact connected oriented
differentiable
surface, $g$ for its genus, and $P$ for a finite nonempty subset of $S$.
Therefore we often suppress $(S,P)$ in the notation and write
$\G$, $A$, $\dots$. We assume that $S-P$ has negative Euler characteristic,
which
amounts to the requirement that if $g=0$, then $|P|\ge 3$.

\head
\section Teichm\"uller spaces
\endhead

\label If $T$ is an oriented $2$-dimensional vector space, then a conformal
structure on $T$ determines an action of the circle group $U(1)$ on $T$ and in
this way $T$ acquires the structure of a $1$-dimensional complex vector space.
Clearly, the converse also holds. Thus, to give the oriented surface $S$ a
conformal structure is  equivalent to give its tangent bundle the structure of
a
complex line bundle. Such a structure comes from a (unique) complex-analytic
structure on $S$, so that $S$ becomes a Riemann surface. By the uniformization
theorem, its universal cover will be isomorphic to the upper half plane.  A
conformal structure on $S$ is given by a section of a fiber bundle whose fibre
is the open convex subset in the vector space of quadratic forms on $\R ^2$
defined by the positive ones. The $C^{\infty}$-topology on this space defines a
topology on the set $\conf (S)$ of conformal structures on $S$. (It also
has a compatible structure of a convex set, so that $\conf (S)$ is
contractible.)

Let $\text{Diffeo} ^+(S\!,\!P)$ denote the group of sense preserving
diffeomorphisms which leave $P$ pointwise fixed, and let
$\text{Diffeo} ^0(S,P)$ denote its identity component. Its ``group of connected
components'',
$$
\G :=\text{Diffeo} ^+(S,P)/\text{Diffeo} ^0(S,P),
$$
is the {\it mapping class group} of $(S,P)$. In this
definition we may replace diffeomorphism by homotopy equivalence (relative $P$)
or all natural choices in between such as PL-homeomorphism,
quasiconformal homeomorphism or plain homeomorphism: we still get the same
group.  Clearly, $\text{Diffeo}^+(S,P)$ acts on the  space of conformal
structures on $S$. The orbit space with respect to its identity component:
$$
\T :=\text{Diffeo} ^0(S,P)\backslash \conf (S)
$$
is called the {\it Teichm\"uller space} of $(S,P)$. It comes with a natural
action of $\G$. If we substitute for $\conf (S)$ the bigger space of conformal
structures inducing the quasiconformal structure underlying the given
differentiable structure and replace $\text{Diffeo}$ by the group of
quasiconformal homeomorphisms of $S$ , then the result is the same. For many
purposes this is actually the most useful characterization.

The Fenchel-Nielsen parametrization shows that $\T $ is
homeomorphic to an open disk. There is even a natural $\G$-invariant
complex-analytic manifold structure on $\T $; if $t\in\T $ is represented by a
Riemann surface $C$ which underlies $S$, the tangent space at $t$ is
identified with $H^1(C,\theta _C(-P))$,  where $\theta _C$ is the sheaf of
holomorphic vector fields on $C$. The action of $\G$ on $\T $ is
properly discontinous and $\G$ has a subgroup of finite index acting
freely (for instance, the kernel of the representation of $\G$ on
$H_1(C;\Z /3)$). This implies that the orbit space
$$
\Mod _g^P:=\G \backslash\T
$$
is in a natural way a normal analytic space with only quotient singularities.

\medskip\label
We can give $\T $ an interpretation as a moduli space: let us first define
an {\it $P$-pointed Riemann surface} $(C,x)$ as a  Riemann surface $C$ together
with an injection $x:P\hookrightarrow C$ such that the automorphism group of
the pair $(C ,x)$ is finite. Say that such an $P$-pointed
Riemann surface  $(C,x)$ is {\it $(S,P)$-marked} if we are given an sense
preserving quasiconformal homeomorphism (henceforth abbreviated as
$\qc$-homeomorphism) $f:S\to C$ which extends $x$, with the understanding
that two such homeomorphisms define the same marking if they are $\qc$-isotopic
relative $P$. Clearly, these markings are permuted in a simply-transitive
manner
by the mapping class group $\G $. An isomorphism of marked
$P$-pointed Riemann surfaces $(C,x,f)$, $(C',x,f')$ is given by an sense
preserving $\qc$-homeomorphism $h:C\to C'$ with $hx=x'$ such that $hf$ is
$\qc$-isotopic to
$f'$ modulo $P$. Now $\T (S,P)$ can be thought of as the space of isomorphism
classes of $(S,P)$-marked Riemann surfaces. So $\Mod _g^P:=\G\backslash\T
(S,P)$ can be identified with the set of isomorphism classes of $P$-pointed
compact Riemann surfaces of genus $g$. It is a coarse moduli space which has a
natural structure of a quasi-projective variety. Knudsen, Deligne and Mumford
showed that there is a distinguished projective completion $\Mod
_g^P\subset\overline{\Mod}{}_g^P$ by the coarse moduli space of stable
$P$-pointed complex curves of genus $g$. (A {\it stable $P$-pointed complex
curve} consists of a complete complex curve $C$ with only simple crossings and
an injection $x$ of $P$ into the nonsingular part of $C$ such that $\Aut (C,x)$
is finite.) It is called the {\it Deligne--Mumford compactification}.

\medskip\label
Let $G=\G /\G _1$ be a finite factor group of $\G $ and put
$$
\Mod _g^P[G]:=\G _1\backslash\T .
$$
Then we have a ramified $G$-covering $\pi _G: \Mod _g^P[G]\to \Mod _g^P$.
The rational cohomology of $\Mod _g^P$ is mapped by $\pi _G^*$ isomorphically
onto the $G$-invariants of the rational cohomology of  $\Mod _g^P[G]$.
If $\G _1$ acts without fixed point on $\T $, then $\T $ can be
regarded as a universal covering space of $\Mod _g^P[G]$, and as $\T $ is
contractible, this implies that $\Mod _g^P[G]$ is a classifying space for $\G
_1$. So the group cohomology of $\G _1$ is the singular cohomology of $\Mod
_g^P[G]$. We get the same statement for $\G $  vis-\`a-vis $\Mod _g^P$,
except that we must use rational coefficients:
$$
H^{\bullet}(\Mod _g^P;\Q )=H^{\bullet}(\G ;\Q ).
$$
This equality represents a gate between algebraic geometry (the left hand side)
and combinatorial group theory (the right hand side).

\head
\section A boundary for Teichm\"uller space
\endhead

We shall give $\T$ a (noncompact) boundary with corners. This is an analogue of
the Borel--Serre compactification for arithmetic
groups and first appeared in a paper by W.J.\ Harvey \cite{\harvey }.

\medskip
We first recall that given a smooth manifold $M$ and a closed submanifold
$N\subset M$ with orientable normal bundle, one has defined the oriented
blowing-up $$ \pi :\Bl _N(M)\to M.
$$
This is a manifold with boundary $\pi ^{-1}N$. The map is an isomorphism
over $M-N$, whereas $\pi ^{-1}N\to N$ can be identified with the sphere bundle
associated to the normal bundle (or more intrinsically, with the bundle of rays
in that bundle) with its obvious projection onto $N$. Notice that in case the
normal bundle has the structure of a complex line bundle, $\pi ^{-1}N\to N$ has
the structure of a $U(1)$-bundle.

This construction generalizes in a straightforward manner to the case  where
$N$
is a union of submanifolds with oriented normal bundles that intersect
multi-transversally; in that case $\Bl _N(M)$ is a manifold with corners and
the
fibres of $\pi$ are products of spheres.

\medskip
Now let $(C,x)$ be a pointed stable curve of genus $g$. Let $\tilde C\to C$ be
its normalization, denote by $\Sigma\subset\tilde C$ the pre-image of
$C_{\sing}$ and consider the composite map
$$
f : \Bl _{\Sigma}(\tilde C)\to\tilde C\to C.
$$
For every $p\in C_{\sing}$, $f^{-1}(p)$ consists of two principal $U(1)$
homogeneous spaces. If we choose for every such $p$ an anti-isomorphism of
these
homegeneous spaces and glue accordingly, then we get an oriented surface over
$C$, $S\to C$, of genus $g$ such that the pre-image of every singular point is
a
circle. We shall interpret the conformal structure on $f^{-1}C_{\reg}$ as a
degenerate conformal structure on $S$.

The choice of the anti-isomorphism over $p$ is the same thing as  the choice of
an anti-isomorphism between $T_pC'$ and $T_pC''$, where $C'$ and $C''$ are the
local branches of $C$ at $p$, given up to a positive real scalar. But this
amounts to choosing a ray in the complex line  $T_pC'\otimes T_pC''$. If we
denote that space of rays by $R_pC$, then our choices are effectively
parametrized by $\prod _{p\in C_{\sing}} R_pC$; this is a principal homogeneous
space of the torus $U(1)^{C_{\sing}}$ that we abbreviate by $R(Z)$.

It is well-kmown that the complex lines $T_pC'\otimes T_pC''$ have an
interpretation in terms of the deformation theory of $C$. Let us recall that
there is a universal deformation
$$
\bigr( (\Cal{C}, C)\to (B,O)\, ;\, x_{\Cal{C}}:(B,O)\times P\to \Cal{C}\bigl)
$$
of $(C,x)$ with as base smooth complex-analytic germ $(B,O)$. Its universal
character implies that the whole situation comes with with an action of the
finite group $\Aut (C,x)$. The $\Aut (C,x)$-orbit space of the base can be
identified with the germ of $\overline{\Mod}_g^P$ at the point defined by
$(C,x)$.

Each singular point $p$ of $C$ determines a smooth divisor $(D_p,O)$ in
$(B,O)$ which parametrizes the deformations of $C$ that do not smooth the
singularity $p$. The fiber over $O$ of the normal bundle of $D_p$,
$T_OB/T_OD_p$, is canonically isomorphic to  $T_pC'\otimes T_pC''$. The
divisors
$D_p$, $p\in C_{\sing}$, intersect with  normal crossings so that their union
$D$  defines an oriented blowing-up: $$ \pi : \Bl _D(B,O)\to (B,O). $$ The
central fiber $\pi ^{-1}O$ is canonically identified with $R(C)$. So over it we
have a canonical family of surfaces of genus $g$. It is easily seen that this
true over all of $\Bl _D(B,O)$, so that we get a family  of oriented genus $g$
surfaces
$$
\Cal{S}\to \Bl _D(B,O).
$$
This family is $P$-pointed.

\smallskip Let $\hat B\to \Bl _D(B,O)$ be a universal cover. Since
$\Bl _D(B,O)$ has the torus $R(C)$ as a deformation retract,  the covering
group
is naturally isomorphic to the fundamental group  of $U(1)^{C_{\sing}}$, i.e.,
to the free abelian group generated by $C_{\sing}$. It is known that the
fundamental group of $U(1)^{C_{\sing}}$ maps injectively to the mapping class
group of a fiber. So the covering transformations permute these markings
freely.

It also follows that $\hat B$ is contractible.  If $\hat{\Cal{S}}\to\hat B$ is
the pull-back of our family of surfaces, then is possible to mark the fibers
simultaneously by means of trivialization $\hat{\Cal{S}}\to S$ relative the
given pointing.   This defines a map from $\hat{B}-\partial\hat B$ to $\T$.
That map is a homeomorphism onto an open subset of $\T$. Now glue $\hat B$ to
$\T$ by means of this map. This clearly endows $\T$ with a partial boundary
with
corners. This can be done over any neighborhood of the Deligne--Mumford
compactification $\overline{\Mod}{}_g^P$ and the essential uniquess of this
construction ensures that the result is a manifold with corners $\hat\T$ whose
interior is $\T$. By construction, $\hat\T$ comes with a $\G$-action that
extends the given one on $\T$. The construction also shows that $\G$ acts
properly discontinuously on $\hat\T$ and that there is a natural proper map
$\G\backslash\hat\T\to \overline{\Mod}_g^P$  whose fibres are finite quotients
of real tori.

There is also a universal family of genus $g$ surfaces over $\hat\T$. As a set,
$\hat\T$ has the following moduli interpretation. Let us define a {\it stable
conformal structure} on $S$ as being given by a closed one-dimensional
submanifold $L\subset S-P$ and a conformal structure on $S-L$ having the
property that contraction of every connected component of $L$ yields a stable
$P$-pointed curve. The set of stable conformal structures is acted on by
$\text{Diffeo} ^+(S,P)$ and the quotient by $\text{Diffeo} ^0(S,P)$ can be
identified with $\hat\T$. The following proposition is well-known and tells us
when a sequence in $\hat\T $ converges.

\proclaim{\label Proposition}
Let $L\subset S-P$ be a compact one-dimensional submanifold such that every
connected component of $S-(P\cup L)$ has negative Euler characteristic. Let
$(J_n)_{n=1}^{\infty}$ be a family of conformal structures on $S$ with the
property that  $(J_n|S-P)_n$ converges uniformly on compact subsets to a stable
conformal structure $J_{\infty}$ on $S$. If $t_{\infty}$ denotes the
corresponding element of $\hat \T$ and $t_n\in\T$ the image of $J_n$, then
$(t_n)_n$ converges to $t_{\infty}$.
\endproclaim

\label In this paper, the space $\hat\T$ will play an auxiliary r\^ole; we will
be more concerned with a quotient $\overline{\T }$ that is a kind of Stein
factorization of the projection $\hat\T\to \overline{\Mod}{}_g^P$:
$\overline{\T }$ is obtained by collapsing every connected component of a fiber
of the latter map to a point. As these connected components are affine spaces
(and hence noncompact in general), the result will not locally compact. Notice
that $\G$ still acts on $\overline{\T }$, and that the orbit space
$\G\backslash\overline{\T }$ can be identified with
$\overline{\Mod}{}_g^P$. So $\overline{\T }\to \overline{\Mod}{}_g^P$ is a
Galois covering with infinite ramification.

\head
\section Quotients of Deligne--Mumford compactifications
\endhead

We introduce certain quotients of $\Mod _g^P$ that are obtained by
identifying points of the boundary of its Deligne--Mumford
compactification and that arise naturally in a combinatorial setting. One such
quotient plays a prominent r\^ole in Kontsevich's proof of a conjecture of
Witten \cite{\konts}. Let us fix a nonempty subset $Q$ of $P$. If $(C,x)$ is an
$P$-pointed stable curve, then the irreducible components of $C$
which contain a point of $Q$ make up a (not necessarily stable) $Q$-pointed
curve $(C_Q,x|Q)$.  The pairs $(C,x)$ for which  every singular point of $C$
lies
on $C_Q$ define a Zariski open subset $U_Q$ of $\overline{\Mod}{}_g^P$. We
define an equivalence relation $R_Q$ on $U_Q$ as follows: two $P$-pointed
stable
curves $(C,x)$ and $(C',x')$ representing points of $U_Q$ are declared to be
$R_Q$-equivalent if there exists an sense preserving homeomorphism $h:C\to C'$
such that $hx =x'$ and $h$ restricts to an analytic isomorphism of $C_Q$ onto
$C'_Q$ as $Q$-pointed curves. We denote its quotient space by $K_Q{\Mod
}_g^P$.
The equivalence relation $R_Q$ has a natural extension $\overline{R}_Q$ to
$\overline{\Mod}{}_g^P$ which is characterized by the property that if we keep
both $C_Q$ and the singular points of $C$ on $C_Q$ fixed, but allow $C$ to
acquire singularities outside $C_Q$, then we stay in the same equivalence
class.
So $K_Q{\Mod }_g^P$ may be regarded as a quotient of $\overline{\Mod}{}_g^P$.

\proclaim{\label Lemma}
The space $K_Q{\Mod }_g^P$ is compact Hausdorff. It contains  ${\Mod }_g^P$ as
an open-dense subset.
\endproclaim
\demo{Proof}
The last assertion of the lemma is  easy and is stated  for the sake of record
only. The first statement is a little ambiguous since it is not clear whether
we
give $K_Q{\Mod }_g^P$ the topology as  a quotient of $U_Q$ or of
$\overline{\Mod
}_g^P$. A priori, the former could be finer than the latter, but we will show
that they are the same. Now $\overline{\Mod}{}_g^P$ is compact and hence so is
every quotient of it. It is therefore enough for us to verify that
$K_Q{\Mod }_g^P$ is Hausdorff as a quotient of $\overline{\Mod}{}_g^P$. This
will be a consequence of the following property of the compactification
$\overline{\Mod}{}_g^P$.

\smallskip
Let $[(C_n,x_n)]_{n=1}^{\infty}$ be a sequence in $U_Q$ converging to $[(C,x)]$
and suppose that all the terms of this sequence have the same topological
type.
Then the intersection of $C_{n,Q}$ with the union of the other irreducible
components of $C_n$ is a finite subset $Z_n$ of the smooth part of $C_{n,Q}$ of
constant cardinality. Let $Z$ be a fixed finite set of this cardinality and
choose for every $n$ a bijection $z_n:Z\cong Z_n$. Then $(C_{n,Q},x_n|Q\sqcup
z_n)$ is a $(Q\sqcup Z)$-pointed curve, which is easily seen to be stable. If
$h$
denotes the arithmetic genus of $C_{n,Q}$, then after passing to a subsequence,
$[(C_{n,Q},x_n|Q\sqcup z_n)]_n$ will converge in $\overline{\Mod}{}_h^{Q\sqcup
Z}$ to some $[(C^*,y\sqcup z)]$. The property alluded to is that $(C^*_Q,
y)=(C_Q,x|Q)$.

\smallskip
To complete the proof, let $[(C_n,x_n)]_{n=1}^{\infty}$ and
$[(C'_n,x'_n)]_{n=1}^{\infty}$ be sequences in $\overline{\Mod}{}_g^P$
converging to $[(C,x)]$ and $[(C',x')]$ repectively such that terms with the
same index are  $R_Q$-equivalent. We must show that $[(C,x)]$ and $[(C',x')]$
are $\overline{R}_Q$-equivalent. But this is immediate from the above mentioned
property. \enddemo

\medskip\label Here is a simple, but perhaps instructive example.  Let $C$ be a
smooth connected projective curve of genus $g\ge 2$. Then $C\times C$
parametrizes a subvariety of  $\overline{\Mod}{} _g^{\{ 0,1\}}$. A point
of the diagonal, $(p,p)\in C\times C$, represents the union of $C$ and $P^1(\C
)$ with $p\in C$ identified with $\infty\in P^1(\C )$ and $i=0,1$ mapping to
$i\in P^1(\C )$. Taking the image in $K_{\{ 0\}}{\Mod } _g^{\{ 0,1\}}$ means
that we disregard the irreducible component $C$ and retain  $P^1(\C )$ with its
three points. So the composite map $C\times C\to K_Q{\Mod } _g^P$ contracts the
diagonal. As A.J.\ de Jong pointed out to me, this contraction can be obtained
algebraically as the normalization of the image of the difference map from
$C\times C$ to the Jacobian of $C$. The contraction can also be realized by the
line bundle on $C\times C$ that is the pull-back of the canonical sheaf under
the projection $(p_0,p_1)\in C^2\to p_0\in C$ twisted by the diagonal (a
positive tensor power of that bundle is without base points).

\medskip\label
Notice that the $\overline{R}_Q$ gets coarser as $Q$ gets  smaller. In
particular, for $Q\subset Q'$, there is a natural quotient mapping
$K_{Q'}{\Mod }_g^P\to K_Q{\Mod }_g^P$.

In this connection we venture the following

\proclaim{Conjecture 1}
The quotients $K_Q\Mod _g^P$  have the structure of a normal
projective variety and such that the quotient map $\overline{\Mod}{}_g^P\to
K_Q\Mod _g^P$ is a morphism.
\endproclaim

If the conjecture holds, then the natural maps
$K_{Q'}{\Mod }_g^P\to K_Q{\Mod }_g^P$, $Q\subset Q'$, are morphisms, too.

We actually expect the corresponding quotient to arise as the image under a
certain linear system without base points. A special (but basic) case is when
$P=Q$ is a singleton:

\proclaim{Conjecture 2} The relatively dualizing sheaf of the universal
stable curve of genus $g\ge 2$, $\overline{\Mod}{}_g^1\to\overline{\Mod}{}_g$,
is semiample, i.e., a positive tensor power of it has no base points.
\endproclaim

S.\ Wolpert \cite{\wolpert} has shown that the natural metric on this
relatively dualizing sheaf has nonnegative curvature and that this curvature
is nonzero in directions transversal to the $R_1$-equivalence classes.
Using this one can show that under the assumption of conjecture 2, a
a positive power of the relatively dualizing sheaf defines a morphism of which
the fibers are the $R_1$-equivalence classes. So conjecture 2 implies
conjecture 1 for the case when $P=Q$ is a singleton.

\medskip\label
These extensions have Teichm\"uller counterparts: for every nonempty $Q\subset
P$ we have a $\G $-equivariant quotient $K_Q\T $ of $\overline{\T
}$ which contains $\T $ and for $Q\subset Q'$ a quotient mapping
$K_{Q'}\T\to K_Q\T$.

It is useful to have a moduli interpretation for these compactifications. We
first remind the reader that one calls a complex-analytic space {\it weakly
normal} if every continuous complex function on an open subset which is
analytic
outside a divisor is analytic. For curves this means that every singular point
with $k$ branches is like the union of the coordinate-axes of $\C ^k$ at the
origin.

We make two definitions:
A {\it $Q$-minimal $P$-pointed curve of genus $g$}  consists of a
connected weakly normal curve $C$, a map $x:P\to C$, and a function  $\epsilon
: C\to\Z _{\ge 0}$ with finite support (the {\it genus defect function}) such
that  \roster
\item $x|Q$ is injective and its image is contained in $C_{\reg}\setminus
x(P-Q)$ and  meets every connected component of that space.
\item The automorphism group of the triple $(C,x,\epsilon )$ is finite
(equivalently: every connected component of $C_{\reg}\setminus
(x(P)\cup\supp\epsilon )$ has negative Euler characteristic).
\item $g=g(\hat C)+\sum_{z\in C} (\epsilon (z)+r(C,z)-1)$, where
$\hat C$ is the normalization of $C$ and $r(C,z)$ is the
number of branches of $(C,z)$
\endroster

The above conditions imply the existence of a continuous map $f:S\to C$ that
extends $x$ such that the pre-image of a point $z\in C$ is connected
submanifold with boundary of $S$ of genus $\epsilon (z)$ with $r(C,z)$ boundary
components if $\epsilon (z)+r(C,z)>1$ and a singleton else. If we are given
such a map up to isotopy relative $P$, then we say that the $Q$-minimal
$P$-pointed curve is {\it marked} by $(S,P)$.

There is an obvious notion of isomorphism: two
$Q$-minimal $P$-pointed curves $(C,x,\epsilon )$, $(C',x',\epsilon ')$ are
declared isomorphic if there exists an isomorphism
$h: C\to C'$ such that $x'h=x$ and $\epsilon 'h=\epsilon$. In the marked
context
we of course also require that $h$ respects the markings.

\proclaim{\label Lemma}
The isomorphism classes of (marked) $Q$-minimal $P$-pointed curves of genus $g$
are in bijective correspondence with the points of $K_Q\Mod _g^P$
($K_Q\T$).
\endproclaim
\demo{Proof}
We content ourselves with indicating how a
$Q$-minimal curve $(C,x,\epsilon)$ determines an element of
$K_Q\Mod _g^P$.  Extend $x$ to a continuous map $f: S\to C$ as above. Let $L$
be the
boundary of $f^{-1}(C_{\sing}\cup\supp\epsilon )$. Now collapse to a point
every
component of $L$ as well as every component of
$f^{-1}(C_{\sing}\cup\supp\epsilon
)$ that is homeomorphic to a cylinder and does not intersect $P$. Then $(\bar
S,\pi x)$ is a stable $P$-pointed pseudosurface. The map $f$ factors through a
map $\bar f:\bar S\to C$ and the irreducible components of $\bar S$ that are
not contracted receive in this way a weakly normal complex structure. Extend
this to a weakly normal complex structure (compatible with the given
orientation) on $\bar S$. Then we get a stable $P$-pointed curve $C$. Its image
in $K_Q{\Mod }_g^P$ only depends on $(C,x,\epsilon)$.
\enddemo

We can form the simplicial scheme $K_{\bullet}{\Mod }_g^P$. Its geometric
realization is a quotient of $\overline{\Mod}{}_g^P$  such that the quotient
map
followed by the structure  map $|K_{\bullet}{\Mod }_g^P |\to\De _P$  is the
projection. We shall show
that $|K_{\bullet}\Mod _g^P|$ is homeomorphic (over $\Delta _P$) to the
semisimplicial complex $\G \backslash A$ that was defined in the
introduction. We look at this complex in more detail in the next section.

\head
\section The arc complex
\endhead

\medskip\label
We consider embedded unoriented loops and arcs $\alpha$ in $S$
which connect two (possibly identical) points of $P$ and avoid all other
points of $P$. In case of a loop we also require that
it be nontrivial in the sense that it does not bound an embedded disk in
$S-P$. Let $\cA $ denote the set of
isotopy relative $P$ of these arcs and loops. We endow this set with the
structure  of an abstract simplicial complex by stipulating that an
$(l+1)$-element subset of $\cA $ defines an $l$-simplex if it is
representable by arcs and loops that do not meet outside $P$. We denote the
geometric realization of this complex by $A$. There is a piecewise linear
map $\lambda$ from $A$ to the simplex $\De _P$ spanned by $P$ characterized
by the property that it sends a vertex $\la\alpha\ra\in \cA $ to the
barycenter of the end points of $\alpha$.  So if $Q$ is a nonempty subset of
$P$
and $\De _Q\subset\De _P$ the corresponding face, then $\lambda^{-1}\De _Q$ is
a
subcomplex of $A$ of which the $0$-simplices may be interpreted as the
isotopy classes of embedded arcs and loops in $S-(P-Q)$ with end points in
$Q$.

\smallskip
We say that the simplex $\langle\alpha _0,\dots ,\alpha _l\rangle$ is {\it
proper} if its star is finite, that is, if it is contained in a finitely many
simplices. This comes down to requiring that each connected component of
$S-\cup
_{\lambda}\alpha _{\lambda}$ is an embedded open disk which contains at most
one
point of $P$. The improper simplices make
up a subcomplex $A_{\infty}\subset A$. We shall denote its complement
$A-A_{\infty}$ by $A^{\circ}$. It is clear that $A$ has an
action of $\G $ which preserves both $A_{\infty}$ and $\lambda$.

\proclaim{\label Lemma}
The group $\Gamma$ has only a finite number of orbits in the set of
simplices of $\cA $. The dimension of a proper simplex is at least
$2g-2+|P|$ and the dimension of every fiber of $\lambda$ is $6g-6+2|P|$.
\endproclaim
\demo{Proof}
The first assertion is a consequence of the fact that
up  to homeomorphism there are only finitely many compact surfaces with an
Euler
characteristic bounded from below (the details are left to the reader).

Let $a=\la\alpha _0,\dots,\alpha _l\ra $ be an $l$-simplex of $A$ and let
$Q\subset P$ the set of points of $P$ that are end point of some $\alpha
_{\lambda}$. This means that $l$ maps the relative interior of $a$ in the
relative interior of $\De _Q$. If $a$ is a proper simplex, then the formula for
the Euler characteristic gives
$$
2-2g=|Q|-(l+1)+d,
$$
where $d$ is the number of connected components of $S-\cup _{\lambda}\alpha
_{\lambda}$. Since every connected component contains at most one point of
$P-Q$,
we have $d\ge |P|-|Q|$. It follows that $l\ge 2g-2+|P|$. If $a$ is maximal in
the pre-image of $\De _Q$, then every connected component of  $S-\cup
_{\lambda}\alpha _{\lambda}$ either is an open disk that contains precisely one
point of $P-Q$ and is bounded by a single member of $a$ or contains no point of
$P-Q$  and is bounded by  three members of $a$. A straightforward computation
shows that then $d={2\over 3}(l+1+|P|-|Q|)$. Substituting this in the formula
for the Euler characteristic gives $l=6g-7+3|Q| +2(|P|-|Q|)=6g-6+2|P|+\dim \De
_Q$.
\enddemo

\remark{Example} We take for $S$ the torus $\R ^2/\Z ^2$ and for $P$ the
origin. An element of $\cA$ is uniquely represented by a circle which is
also a subgroup of $S$. Such a subgroup is the image of a line in $\R ^2$
through  the origin and another point of $\Z ^2$. In this way we obtain
an identification of $\cA$ with the rational projective line
$\bold{P}^1(\Q)$. The two circles defined by the relatively prime pairs of
integers $(x_0,x_1)$ and $(y_0,y_1)$ define a $1$-simplex iff they do not
intersect outside the origin. This is the case iff
$x_0y_1-x_1y_0\not=\pm 1$, or equivalently, iff $x=(x_0,x_1)$ and $y=(y_0,y_1)$
make up a basis of $\Z ^2$. Then this $1$-simplex is adjacent to exactly two
$2$-simplices, namely those defined by $\{ x,y,x+y\}$ and $\{ x,y,x-y\}$.
A simplex is proper iff
it is of dimension $>0$. The geometric realization of $A$ can be pictured in
the upper half plane (with the vertex at $\infty$ missing) as a hyperbolic
tesselation associated to a subgroup of the modular group of index two.
\endremark

\midspace{50mm}\caption{\fig The arc complex of a once-pointed torus}

Let $b\cA $ denote the barycentric subdivision of $\cA $. So a
vertex of $b\cA $ is the barycenter of a simplex $a$ of $\cA $ and a
$k$-simplex of $b\Cal{A}$ is spanned by the barycenters of a strictly
increasing chain $a_0<a_1<\cdots <a_k$ of simplices of $b\Cal{A}$. Let
$\Cal{A}_{\pr}$ denote the full subcomplex of $b\Cal{A}$ whose
vertices are the barycenters of proper simplices. Clearly, its geometric
realization $A_{\pr}$ can be viewed as a subset of $A^{\circ}$. In the
previous example we have drawn $A_{\pr}$ with dotted lines.

\proclaim{\label Proposition}
The fibres of $\lambda |A_{\pr}$ have dimension $4g-4+|P|$ and
there is a natural $\G $-equivariant deformation retraction of
$A^{\circ}$ resp. $A-A_{\pr}$ onto $A_{\pr}$ resp.
$A_{\infty }$ which preserves the pre-image of every relatively open
face of $\De _P$ under $\lambda$.  \endproclaim
\demo{Proof}
A $k$-simplex of $\Cal{A}_{\pr}$ is represented by a chain
$a_0<a_1<\cdots <a_k$ of simplices of $\Cal{A}$ with $a_0$ proper.
According to the previous lemma $\dim a_0\ge 2g-2+|P|$ and $\dim a_k\le
6g-6+2|P|+\dim\De _Q$, where $Q\subset P$ is the smallest subset of $P$ such
that $\lambda$ maps $a_k$ in $\De _Q$. So $k\le (6g-6+2|P|+\dim \De
_Q)-(2g-2+|P|)=4g-4+|P|+\dim \De _Q$. It is easily verified that this value is
attained.

The proof of the remaining assertions is a standard argument in the theory of
simplicial complexes, but let us give it nevertheless, say for
$A_{\pr}\subset A^{\circ}$. If  $x\in
A^{\circ}=bA-bA_{\infty}$, then we can write $x=\sum
_{i=0}^k x_ia_i$ with $a_0<a_1<\cdots <a_k$, $x_i>0$, and $a_k$ proper. Let $r$
be the first index such that $a_r$ is proper. Then
$$
x':=\sum _{i=r}^k (\sum _{j=r}^k x_j)^{-1}x_ia_i\in A_{\pr}
$$
and $x(t):=(1-t)x+tx'$ defines a deformation retraction of
$A^{\circ}$ onto $A_{\pr}$.
\enddemo

Our goal is to construct a $\G $-equivariant homeomorphism of
$A$ onto $|K_{\bullet}\T |$ which commutes with the given projections
onto $\De _P$. For this we first  need to discuss ribbon graphs.

\head
\section Ribbon graphs
\endhead

\label A {\it ribbon graph} is a nonempty finite graph in which we allow
loops and multiple bonds, but not isolated points (in other words, a
semi-simplicial complex of pure dimension $1$), such that for every
vertex we are given a cyclic order of its outgoing edges.

A finite graph embedded in an oriented surface acquires
naturally such a structure. Conversely, a ribbon graph can be embedded in an
oriented surface of which it is a deformation retract. For instance,

\midspace{50mm}\caption{\fig Ambient surface of a ribbon graph}

\noindent This surface can be compactified by adding a finite number of points
so that the result is a surface.

This compactification can be obtained in a purely
combinatorial way as follows.
Let $G$ be a ribbon graph. Denote by $X(G)$ its set of oriented
edges (so that each edge determines two distinct elements of $X(G)$).
Reversal of orientation defines a fixed point free involution $\sigma _1$ in
$X(G)$. For $e\in X(G )$, let $v$ be its vertex of origin, and
denote by $\sigma _0(e)\in X(G )$ the outgoing edge of $v$ that succedes
$e$ relative the given cyclic order. This defines a permutation $\sigma _0$ of
$X(G )$. We define the permutation $\sigma _{\infty}$ by the equality
$\sigma _{\infty}\sigma _1\sigma _0=1$.

\midspace{50mm}\caption{\fig The operations $\sigma _i$}

Denote the orbit space of
$\sigma _i$ in $X(G )$ by $X_i(G )$. For $i=0$ resp. $i=1$ it can
be identified with the set of vertices resp. of (unoriented) edges of $G$;
the elements of  $X_{\infty}(G )$ are called {\it boundary cycles}. So
$G $ can be reconstructed from $X(G )$ equipped with the permutations
$\sigma _0$ and $\sigma _1$. (Indeed, any nonempty finite set equipped with a
fixed point free involution and another permutation determines a ribbon
graph.)

\medskip\label
Let $K$ be the two-simplex with vertices $v_0,\bar v_0,v_{\infty}$ with
the orientation given by this order. The midpoint of the face $\la v_0,
\bar v_0\ra$ is denoted $v_1$. Denote by a ``bar'' the involution of $K$
which interchanges $v_0$ and $\bar v_0$ and leaves $v_{\infty}$ (hence also
$v_1$) fixed. We define a semi-simplicial complex $S(G )$ as a quotient of
$K\times X(G )$ by identifying the oriented $1$-simplices $\la v_0,\bar
v_0\ra \times \{e \}$ with  $\la\bar v_0,v_0\ra\times\{ \sigma _1e\}$ and $\la
v_0,v_{\infty}\ra \times\{ e\}$ with $\la \bar v_0,v_{\infty}\ra \times\{
\sigma_0 e\}$.

Since the the disjoint union of the $X_0(G)$ and
$X_{\infty}(G )$ appears here as the set of $0$-simplices, we will often regard
these two as subsets of $S(G)$. In what follows a special r\^ole is played by
the $0$-simplices that either belong to $X_{\infty}(G)$ or are a vertex of $G$
of valency $\le 2$. We shall call such points {\it distinguished}.

We shall write $K_e$ for the image of $K\times
\{ e\}$ and we call it the {\it tile} defined by $e$. The full subcomplex
spanned by $X_0(G )$ can be identified with $G$, see the picture
below.

\midspace{50mm}\caption{\fig Combinatorial construction of the ambient surface}

It is not difficult to see that the geometric realization of $S(G )$ is a
compact surface. The given orientation of $K$ determines one  of $S(G )$
and this orientation is compatible with the ribbon graph structure of $X(G
)$. The surface has a piecewise linear structure and hence a quasiconformal
structure.

The image of $\la v_1,v_{\infty}\ra\times X(G )$ is the barycentric
subdivision of another ribbon graph, called the dual of $G$, and denoted by
$G ^*$. (It is essentially obtained by passing from $(X(G );\sigma
_0,\sigma _1)$ to $(X(G );\sigma _{\infty},\sigma _1)$ and using a natural
identification of $S(G ^*)$ with $S(G )$.)  Observe that the edges of
$G ^*$ are indexed by the edges of $G$.

{\it Remark.} The permutations $\sigma _0,\sigma _1,\sigma _{\infty}$
associated
to a ribbon graph $\G$ arise as monodromies in the following manner. Let $S_0$
be the topological sphere obtained from $K$ by identifying points on its
boundary according to the  involution ``bar'' and denote the image of $v_z$ by
$z\in S_0$ ($z=0,1,\infty$). It is clear that there is a natural finite
quotient
map $S(G )\to S_0$. This map is ramified covering which branch locus $\{
0,1,\infty\}$.  The restriction to $K^{\circ}\hookrightarrow S_0$ is naturally
identified with $K^{\circ}\times X(G )$ and the monodromy of $S(G )\to S_0$
around $z\in \{0,1,\infty\}$ is given by the permutation $\sigma _z$ acting on
the second factor.

\head
\section Metrized ribbon graphs
\endhead

\medskip\label
A {\it metric} on a ribbon graph is $G $ simply a map from its edges to $\R
_{>0}$.  If this map has in addition the property that the total length of the
graph is $1$, then we call it a {\it unital metric}.

A {\it conformal structure} on $G $ is a metric on every connected component of
$G $, given up to a factor of proportionality. This is of course equivalent to
be given a unital metric on every connected component of $G $. We denote the
space of conformal structures on $G $ by $\conf (G  )$. So for connected $G $,
$\conf (G
)$ may be identified with the open simplex spanned by the set of edges of $G
$.

\medskip\label
Let $r:K\to [0,1]$ be  the barycentric coordinate which is $1$ in
$v_{\infty}\in
K$ and $0$ in $v_0$ and $\bar v_0$ and identify $K-\{\infty \}$ with $\la
v_0,\bar v_0\ra \times\R _{\ge 0}$, where the first component is an obvious
projection and the second is given by $-\log r$.  Suppose that we are given a
ribbon graph $G $ with metric $l:X_1(G )\to\R _{>0}$. This determines a
complete
piecewise Euclidean metric on $S(G  ) - X_{\infty }(G  )$ as follows: give
$(K-\{\infty \})\times \{e\}$ the metric which under its identification with
$\la v_0,\bar v_0\rangle \times\R _{\ge 0}$ corresponds to the translation
invariant product metric for which $\la v_0,\bar v_0\ra $ has length $l(e)$ and
the second component has the standard metric. This descends to a metric on $S(G
) - X_{\infty }(G  )$. The complement of the vertex set of $S(G  )$
has a unique smooth structure for which this metric is Riemannian on
that set. It is easy to check that its underlying conformal structure extends
across the vertices, so that now $S(G  )$ acquires a conformal structure.
We denote the Riemann surface thus obtained by $C(G  ,l)$. This Riemann surface
comes with a meromomorphic quadratic differential $q_l$ whose absolute value
gives the metric: if we identify the interior of the tile $K_e$ as a metric
space in the obvious way with the Euclidean rectangle $\{ z\in \C: \Im (z)>0,
|\Re (z)|< {1\over 2}l(e)\}$, then this is a complex-analytic chart and the
quadratic differential is given by $dz\otimes dz$. One finds that $q_l$ has a
pole of order two at each point of $X_{\infty}(G  )$ and a zero at of order
$k-2$ at each $k$-valent vertex of $G $ (so a pole of order one at a univalent
vertex). This implies that successive outgoing oriented edges at a $k$-valent
vertex make an angle of $2\pi/k$. There are no other singularities of $q_l$.
Observe that as a piecewise-linear complex valued quadratic differential on
$S(G  )$, $q_l$ embodies all the extra structure: the smooth structure, the
metric and (hence) the complex-analytic structure.

Notice that the conformal structure on $S(G  )$ only depends on the conformal
structure on $G $ subordinate to $l$. Hence we can always assume that $l$ is
unital on every connected component of $G $.

If $v$ is a bivalent vertex of $G $, then ``forgetting'' that vertex yields a
metrized ribbon graph of which the associated Riemann surface can be identified
with $C(G  ,l)$.

\smallskip
In case of a $l$ is a metric on a partial ribbon graph we do essentially the
same construction where the metric on the incomplete edges should be thought of
as having the value $\infty$.  So if $e$ is an oriented edge without end point,
then $K^{\circ}_e$ and $K^{\circ}_{\sigma _1(e)}$ are Euclidean quadrant
(isomorphic to $\R _{\ge 0}^2$) such that $e$ corresponds to the positive
$x$-axis in the former and to the positive $y$-axis in the latter. Again on
verifies that the underlying conformal structure extends across the vertices so
that we find a Riemann surface $C(G  ,l)$. This time the quadratic differential
$q_l$ may have higher order poles at the points of $X_{\infty}(G  )$. In fact,
the pole order at $\beta\in X_{\infty}(G  )$ will be $2$ plus the number of
incomplete (unoriented) edges that occur in $beta$.

\medskip\label
An {\it $P$-pointed ribbon graph} is an ribbon graph $G $ together with an
injection $x: P\hookrightarrow  X_{\infty}(G)\sqcup X_0(G )$ whose image
contains all the distinguished points. Notice that in that case every connected
component of $S(G)-x(P)$ has negative Euler characteristic: this is because
$S(G)-X_{\infty}(G)$ admits $G$ as a deformation retract and every connected
graph which is contractible (resp.\ a homotopy circle) has at least two (resp.\
one) vertices of valency at most $2$.

\medskip
Let $(G,x)$ be an $P$-pointed ribbon graph. If $s$ is an edge of $G $ which
is neither isolated nor a loop, then collapsing that edge yields a ribbon graph
$G  /s$. It inherits an $P$-pointing iff not both of its vertices are in the
image of $P$. The corresponding surface $S(G  /s)$ is obtained as a quotient of
$S(G  )$ by collapsing the two tiles defined by $s$ according to the level sets
of $r$. We call this an {\it edge collapse}.

If $s$ is a non-isolated loop, and for some orientation $e$ of $s$, $e$ is by
itself a boundary cycle, then it is still true that $G  /s$ is a ribbon graph.
In this case, $G  /s$ inherits a $P$-pointing ifand only if the vertex of $s$
is
not in the image of $P$. The surface $S(G  /s)$ is then obtained by collapsing
$K_e$ to a point (a {\it total collapse}) and by applying an edge collapse to
the
opposite tile $K_{\sigma _1e}$.

In either case the quotient map $S(G)\to S(G/s)$ has in its homotopy class
relative $P$ a unique isotopy class relative $P$ of $\qc$-homeomorphisms.

\smallskip
We can apply these two procedures successively to a collection $Z$ of edges of
$G $ if and only if every connected component of the corresponding subgraph $G
_Z\subset G$ is
\roster
\item either a tree with at most one  marked vertex or
\item a homotopy circle without marked vertices which contains an entire
boundary cycle of $G $.
\endroster
We then say that $Z$ is {\it negligible}.
So if $Z$ is negligible and $G /G  _Z$ is the semi-simplicial complex  obtained
by collapsing  every connected component of $G  _Z$ to a point, then $G /G  _Z$
has still the structure of a ribbon graph pointed by $P$ and the corresponding
surface $S(G/G _Z)$ can be obtained by means of a succession of edge
collapses and contractions of the tiles labeled by the oriented edges in $Z$.
The quotient map $S(G)\to S(G/G _Z)$ determines an isotopy class relative
$P$  of sense preserving $\qc$-homeomorphisms $S(G  )\to S(G  /G  _Z)$.

\smallskip
An {\it almost-metric} on $G $ is a function $l: X_1(G  )\to \R _{\ge 0}$ whose
zero set $Z$ is negligible. It is clear that $l$ then factorizes over a
metrized
ribbon graph $G  /G  _Z$ with metric (still denoted) $l$ and we define
$C(G  ,l)$ simply as $C(G  /G  _Z, l)$. We have a corresponding notion of an
{\it almost-conformal structure}.

Denote the space of unital almost-conformal structures on $(G,x)$ by  $\aconf
(G ,x)$. It is clear that for a negligible $Z\subset X_1(G  )$, we have a
natural embedding of $\aconf (G  /G  _Z,x)$ in $\aconf (G,x)$.

\medskip\label We now assume that $G $ is a connected  ribbon graph.
Over $\conf (G  )$ lives a ``tautological'' topologically trivial family of
metrized graphs and a corresponding family of Riemann surfaces. We extend the
latter as a family of pseudosurfaces; in section 8 we give each of its fibers
the structure of a weakly normal curve.

The family appears as a factor of the projection $S(G)\times a(G)\to a(G)$ and
is defined as follows. Any edge $s$
of $G $ determines by definition a vertex  of $a(G)$.  The
codimension-one face opposite this vertex is identified with $a(G/s)$ and for
each orientation $e$ of $s$, we apply an edge collapse to $K_e\times a(G/s)$
relative its projection onto $a(G/s)$. Likewise, every boundary cycle $\beta$
of $G$ determines a face $a(G/G_{\beta})$ of $a(G)$ and we perform a total
collapse on the tiles $K_e\times a(G/G_{\beta})$ relative $a(G/G_{\beta})$ with
$e\in\beta$. The result is a semisimplicial space $\CC (G)$ that comes with a
projection $\pi _G:\CC (G)\to a(G)$.

Over $l\in \conf (G  )$ the fiber is the surface $S(G)$; it has a conformal
structure which makes it canonically isomorphic to $C(G,l)$. That last fact is
still true in case $l\in \aconf (G)$. The fiber $\CC (G)_l$ over an arbitrary
$l\in a(G)$ is gotten as follows.  Let $Z\subset X_1(G  )$ be the zero set of
$l$ and let  $S(G  )_Z$ be the quotient of $S(G  )$ obtained by performing for
every oriented edge $e$ of $Z$ a contraction or an edge collapse on $K_e$,
depending on whether or not the boundary cycle of $G$ generated by $e$ is
contained in $G_Z$.  Then $\CC (G)_l$ can be identified with $S(G)_Z$. We
will see in section 8  that $S(G)_Z$ is a pseudosurface and that $\CC (G)_l$
has a natural conformal structure on its smooth part given by quadratic
differential. (This conformal structure determines a unique complex-analytic
structure such that $\CC (G)_l$ is weakly normal.)

\medskip\label
We conclude this discussion with a few remarks.

Every element  of $X_0(G)\sqcup X_{\infty}(G)$ determines a section
of $\CC (G)\to a(G)$. Those that are indexed by $P$ are disjoint over
$\aconf (G)$.

One can show that the complement of the sections defined by the elements of
$X_0(G)\sqcup X_{\infty}(G)$ has a natural smooth structure. (To see this, use
an atlas naturally indexed by the elements of $X_1(G)\sqcup X_{\infty}(G)$.)
The conformal structures along the the fibers vary differentiably on this open
subset.

\head
\section Moduli spaces
\endhead

\medskip\label
We say that a ribbon graph $G $ is {\it $(S,P)$-marked}
(or briefly, {\it marked}) if we are given a given
isotopy class relative $P$ of sense preserving $\qc$-homeomorphisms $f: S\cong
S(G)$ such that $f|P$ defines a $P$-pointing of $G $: $f$ maps $P$ to
$X_{\infty}(G)\sqcup X_0(G)$ and its image contains the distinguished points.
It is clear that $G$ permutes the markings.

\medskip\label
We claim that a marked ribbon graph is the same thing as a proper
simplex of $\Cal{A}$. Let $f :S\cong S(G)$ be a marking. Regard the dual
ribbon graph $G  ^*$ as lying on $S(G  )$. Then the pre-image of every edge of
$G^*$ under $f$ connects two points of $P$ and therefore the collection of
these determines a simplex $a(G,f)$ of $\Cal{A}$. A connected component
of $S-G ^*$ is given by a vertex of $G $; it contains one or no point of $P$
depending on whether this vertex is marked by $P$. If the vertex is unmarked it
has valency $k\ge 3$ and the connected component is $k$-gon.
So distinct edges of $G^*$ yield distinct vertices of $a(G,f)$ and
$a(G,f)$ is a proper simplex. We also notice that the space of unital metrics
$\conf (G  )$ may be identified with the relative interior of $a(G,f)$; we
shall
therefore denote that relative interior by $\conf (G,f)$.

Conversely, if $a=\la\alpha _0,\dots ,\alpha _l\ra$ is a proper simplex of
$A$, then the union of the $\alpha _i$'s define a ribbon graph $G  _a$ on
$S$ with vertex set contained in $P$. It is easily seen that the inclusion $G
_a\subset S$ extends to a $\qc$-homeomorphism $S(G_a)\to S$ such that
$X_{\infty}(G_a)$ is mapped in $P$. If we identify $S(G_a ^*)$ with
$S(G  ,S)$, then we see that $G_a$ has in a natural way the structure of a
marked ribbon graph.

We remark that $\conf (G,f)$ has maximal dimension iff all vertices
of $G$ are trivalent (so that $P$ maps bijectively onto the set boundary cycles
of $G$).

\proclaim{\label Lemma}
Let $a$ be a proper simplex of $A$ as above with associated
marked ribbon graph $(G  ,f)$. Let $Z\subset X_1(G  )$ be a set of
edges of $G $ and let $a(G /G  _Z)$ be the codimension $|Z|$ face of $a$
opposite the face defined by $Z$. Then $Z$ is negligible if and only if
$a(G /G  _Z)$ is proper and in that case $S(G  /G  _Z)$
inherits an marking (denoted $f/Z$).
\endproclaim
\demo{Proof} It is enough to show this in case $Z$
has only one element and this we leave to the reader.
\enddemo

So given a marking $f$, then the space of unital almost-metrics
$\aconf (G  ,f|P)$ may be identified with  $|a(G  ,f)|\cap A^{\circ}$. We
denote the latter by $\aconf (G  ,f)$.

The restriction of $\lambda :A\to\De _P$ to $\aconf (G  ,f)$ has the following
simple description: for $p\in P$ the corresponding barycentric
coordinate $\lambda _p$ is in case $f(p)$ corresponds to a boundary cycle, half
the length of that cycle and it is zero otherwise.

\smallskip
Remember that
every proper simplex of $A$ is of the form $a(G,f)$ and that over such a
simplex we have defined in section 6 the family $\CC (G  )\to a(G  ,f)$. As
each inclusion of proper simplices is canonically covered by an inclusion of
the
corresponding families, this gives us a global family $\pi :\CC \to
A$. This family comes with sections labeled by $P$.

Summing up:

\proclaim{\label Proposition}
The set of points of $A^{\circ}$ is naturally interpreted as
the set isomorphism classes of marked ribbon graphs endowed with a
unital metric. It is obtained from the spaces $\aconf (G  ,f)$ by identifying
$\aconf (G
/G  _Z,f/Z)$ with its image in $\aconf (G  ,f)$ for every negligible $Z\subset
X_1(G
)$. Moreover, $A$  supports a family $\pi :\CC \to A$ of weakly
normal curves with sections indexed by $P$. Over $A^{\circ}$ these sections
are disjoint, the family is locally trivial with fiber $S$ and each fiber comes
with a complex structure which varies continuously with the base point.
\endproclaim

In the next section  we shall discuss the fibers over $A_{\infty}$.\par

The family $\pi $ restricted to $A^{\circ}$ defines a classifying map
$\Phi : A^{\circ}\to\T $. This map is continuous and clearly
$\G$-equivariant. The following theorem is a rather direct consequence of
the work of Strebel.

\proclaim{\label Theorem}
The map
$$
\Psi ^{\circ}:=(\Phi,\lambda) :A^{\circ}\to\T \times\De _P
$$
is a homeomorphism.
\endproclaim

The observation that Strebel's work leads to theorems of this type is due to
Thurston, Mumford and Harer \cite{\harera}. (We did not come across this
version, though.)

For the proof we must discuss Jenkins-Strebel differentials first. Let  $R$ be
a
Riemann surface. If $q$ is a meromorphic quadratic differential on $R$, then at
each point $p$ of $R$ where $q$ has neither a zero nor a pole the tangent
vectors at $p$ on which $q$ takes a real value $\ge 0$ form a real
line in $T_zC$. This defines a foliation on $R$ minus the singular set of $q$.
If the union of the closed leaves of this foliation is dense in $R$,
then $q$ is called a {\it Jenkins-Strebel differential}. Suppose $q$ is such a
differential. Then a local consideration shows that $q$ has no poles of order
$>2$ and that the double residue at a pole of order $2$ is a negative real
number. The form $q$ determines a Riemann metric $|q|$ on the complement  of
the
singular set of $q$.  This metric is locally like $|dz|^2$ and hence flat.  The
union $K$ of the non-closed leaves and the singular points of $q$ of order $\ge
-1$ is closed in $R$. It is an embedded graph with a singularity of order $k$
being a vertex of valency $k+2$; it is called the {\it critical graph} of $q$.
Each connected component of the complement of $K$ is either a flat annulus
(metrically a flat cylinder) or a disk containing a unique pole of order two
(metrically outside this pole a flat semi-infinite cylinder) or a copy of $\C
-\{ 0\}$.

Suppose that $R$ is the complement of a finite subset of a compact  Riemann
surface $C$. Then $q$ is also a  Jenkins-Strebel differential on $C$ and the
closure $\overline{K}$ of $K$ in $C$ is an embedded graph. (When $C$ has genus
zero it may happen that this closure becomes a closed orbit on $C$, so
$\overline{K}$ may depend on $R$. It can be shown however, that this is the
only
such case.) Clearly, $\overline{K}$ has the structure of a ribbon graph. Notice
that $q$ defines a metric on it.

\proclaim{\label Theorem} {\rm(Strebel)}
Let $(C,x)$ be a compact connected $P$-pointed Riemann surface such that
is not the two-pointed Riemann sphere and let $\lambda\in\Delta _P$. Then there
exists a Jenkins--Strebel differential $q$ on $C$ with the property
that the union of the closed leaves of $q$ form semi-infinite
cylinders around the points of $x(p)$ with $\lambda (p)\not= 0$ (of
circumference $\lambda (p)$) and the points $x(p)$ with $\lambda (p)=0$ lie on
the critical graph of $q$. Moreover, such a $q$ is unique.
\endproclaim
\demo{Proof}
Denote by $Q\subset P$ denote the zero set of $\lambda$ and put
 $Q':=P-Q$. If $|Q'|\ge 2$, then the asserted properties follow from Theorem
$23.5$ of \cite{\strebel} applied to the Riemann surface $C-x(Q)$ with
circumferences given by $p\in Q'\mapsto \lambda (p)$. (The fact that $q$ will
have at the points of $Q$ order $\ge -1$ follows from the discussion above.) In
case $Q'$ is a singleton $\{ p\}$, then Theorem $23.2$ of \cite{\strebel}
implies that there is Jenkins--Strebel differential on $C-x(Q)$ for which all
the closed leaves belong to the cylinder about $p$. This differential is unique
up to a positive real scalar factor and hence the theorem follows in this case,
too.
\enddemo

We shall refer to $\lambda$ as a {\it circumference function} of
 $(C,x)$, the name being suggested by the above theorem. So such a function
determines a metrized ribbon graph $(G _{\lambda},l_{\lambda})$ in $C$
(denoted by $\overline{K}$ in the discussion above). Notice that if $\lambda
(p)=0$, then $x(p)$ is a univalent vertex  or an interior
point of an edge of $G _{\lambda}$; if $\lambda (p)\not=  0$, then $x(p)$
defines a boundary cycle of $G _{\lambda}$. Moreover, all univalent vertices
and boundary cycles of $G _{\lambda}$ are thus obtained. In other words, $G
_{\lambda}$ is in a natural manner an $P$-pointed ribbon graph. The associated
$P$-pointed curve $C(G _{\lambda},l_{\lambda})$ is canonically isomorphic to
$(C,x)$: this is clear on the complement of the union of $x(P)$ and the vertex
set of  $G _{\lambda}$. Hence it is true everywhere.

\demo{Proof of \refer{7.5}}
The above discussion shows that $\Psi ^{\circ}$ has a unique
inverse, in other words, that it is bijective. Since $\Psi ^{\circ}$ is
continuous and has  locally compact domain and range, it must be a
homeomorphism.
\enddemo

\proclaim{\label Corollary} {\rm (Harer \cite{\harera})}
For nonempty $P$, the moduli space $\Mod _g^P$ has
 the homotopy type of a finite semi-simplicial complex of dimension $\le
4g-4+|P|$. In particular,  $\Mod _g^P$ has no homology or cohomology in
dimension $>4g-4+|P|$.
\endproclaim
\demo{Proof}
Choose $p\in P$ and regard $p$ as a vertex of $\De _P$.  Then $\T $
is by \refer{7.5} equivariantly homeomorphic to  $\lambda ^{-1}(p)\cap
A^{\circ}$. Now apply \refer{4.3}.
\enddemo

\head
\section Minimal models
\endhead

In this section we introduce a combinatorial analogue of a $Q$-minimal
$P$-pointed curve. Here $(G ,x)$ is a connected marked ribbon graph.

\medskip\label
We say that a set $Z$ of edges of $G$ is {\it semistable} if no component of $G
_Z$ is the set of edges of a negligible subset and every univalent vertex of
$G_Z$ is in the image of $x$. Then every component of $G_Z$ which  is
contractible contains at least two vertices in $x(P)$. A component which is a
homotopy circle without a vertex in $x(P)$ is necessarily a
topological circle which is not a boundary cycle of $G$. It
is clear that every subset $Z\subset X_1(G)$ has a maximal semistable subset
$Z^{\sst}$. Notice that $Z^{\sst}-Z$ is a negligible subset of $X_1(G)$ so that
if we put $G':=G/G_{Z^{\sst}-Z}$, then $S(G')$ is $\qc$-homeomorphic relative
$P$ to $S(G)$. We sometimes regard $G_{Z^{\sst}}$ as a graph on $S(G')$, so
that with this convention $G/G_Z=G'/G'_{Z^{\sst}}$.

\medskip\label Let be given a proper subset $Z$ of $X_1(G)$.
We can associate to $Z$ two ribbon graphs: one with edge set $Z$ and another
with
edge set $X_1(G)-Z$. In the first case we give $G _Z$ an induced structure of
ribbon graph by telling how the corresponding operator $\sigma _0$  acts on
$X(G _Z)$: it sends $e\in X(G _Z)$ to the first term of the sequence  $(\sigma
_0^k(e))_{k\ge 1}$ which is in $X(G _Z)$.  The second case is in a sense dual
to the first: we define a ribbon graph $G/G _Z$ with $X_1(G )-Z$ as its set
of edges and the corresponding operator $\sigma _{\infty}$ sends $e\in X(G
)-X(G _Z)$ to the first term  of the sequence $(\sigma _{\infty}^k(e))_{k\ge
1}$ which is not in $X(G _Z)$. This ribbon graph naturally maps onto a subgraph
of $G$, but this map need not be injective as it may identify distinct vertices
of $G/G_Z$.

A vertex of $G/G_Z$ that is in the image of an oriented edge in $Z^{\sst}$ will
be called {\it exceptional}. Any such vertex corresponds to a boundary cycle of
$G_{Z^{\sst}}$ that is not a boundary cycle of $G$ (and vice versa), reason for
us to call such boundary cycles {\it exceptional} also.

\proclaim{\label Lemma}
There is a natural identification mapping of $S(G/G_Z)\to S(G)_Z$. This
map identifies two distinct points if and only if both are exceptional vertices
of $S(G/G _Z)$ that come from a boundary cycle of the same component of $G
_{Z^{\sst}}$. In particular, $S(G)_Z$ is a pseudosurface whose combinatorial
normalization is $S(G/G_Z)$. Moreover, every distinguished point of $G/G_Z$
comes from a distinguished point of $G$ or is exceptional.

\endproclaim
\demo{Proof}
Straightforward.
\enddemo

In this situation we have a genus defect function $\epsilon : S(G)_Z\to \Z
_{\ge 0}$ which assigns to the image of an exceptional vertex the genus of the
corresponding component of $S(G_{Z^{\sst}})$ and is zero else.

\medskip\label
Choose an $l\in a(G)$.
In \refer{6.4} we constructed a map $\pi _{G }:\CC (G )\to a(G )$ and we
noticed that that the fiber over $l$, $\CC (G )_l$, can
be identified with $S(G )_Z$, where $Z$ is the zero set of $l$. Since $l$
determines a unital metric on $G/G_Z$, we have a Riemann surface
$C(G/G_Z,l)$ with underlying space $S(G/G_Z)$. We use the previous lemma
to give $\CC (G)_l$ the unique complex-analytic structure for which $\CC (G
)_l$ is weakly normal and $C(G/G_Z,l)\to C(G)_l$ is its normalization.

\proclaim{\label Proposition}
Let $Q$ be the set of $p\in P$ that map to a boundary cycle of $G$ of positive
length. Then $(Q,\epsilon ,P\to S(G )\to \CC (G )_l)$ give $\CC (G )_l$ the
structure of a $Q$-minimal $P$-pointed curve.

\endproclaim
\demo{Proof}
We verify the defining properties of \refer{3.4}. The property for $p\in P$ to
belong to $Q$ is equivalent to $x(p)\in X_{\infty}(G /G _Z)$. The first
property now follows.
For the second we must show that $S(G /G _Z)-X_{\infty}(G /G
_Z)-\{\text{exceptional vertices}\}$ has negative Euler characteristic. But
this follows from the fact that this is (by \refer{8.3}) just the complement of
the set of distinguished points on $G/G_Z$. The verification of the third
property is left to the reader.
\enddemo

Suppose we are given a marking $f$ of $G$ that extends the pointing by $x$.
This determines a marking of $\CC (G )_l$ by $(S,P)$.
In view of the moduli interpretation \refer{3.5}, the structure present on
$\CC (G )_l$ determines a point of $K_Q\T$. By letting $l$ vary over
the elements of $a(G,f)$, we thus obtain a map $a(G,f)\to |K_{\bullet}\T |$
commuting with the given  maps of domain and range to $\Delta _P$. For a
negligible edge $s$ of $G$ the restriction of this map to $a(G /s,f/s)$
coincides with the one defined for that simplex. This results in an
$\G$-equivariant map  $\Psi :A\to |K_{\bullet}\T |$. We can now state our
first main result. It gives an analytic interpretation of $A$:

\proclaim{\label Theorem}
The map $\Psi : A\to |K_{\bullet}\T |$ is a $\G$-equivariant continuous
bijection that commutes with the given maps to $\De _P$.
\endproclaim

The main difficulty is to show that $\Psi$ is continuous.
We postpone the proof to a point where we have treated the combinatorial
version of the Deligne--Mumford compactification.
The reader may wonder whether $\Psi $ is a homeomorphism. The answer is that it
is not, as is illustated by the case $g=1$, $P$ a singleton: then
$|K_{\bullet}\T |$ is the union of the
upper half plane and $P^1(\Q )$. Near $\infty$ it has the horocyclic topology
but the topology it receives from its triangulation is much
finer: a subset of the upper half plane is the complement of a neighborhood of
$\infty$ if and only if its intersection with any vertical strip of bounded
width is bounded.

\head
\section Stable ribbon graphs
\endhead

Here we introduce the ribbon graph analogue of a stable
$P$-pointed curve. That our definition is the natural one may not be
immediately obvious, but that this is indeed the case will become apparent in
the discussion following the definition and in section 10.

\medskip\label
Suppose we are given a ribbon graph $G$ and an injection $x:P\to X_0(G)\sqcup
X_{\infty}(G)$. We no longer assume that $x(P)$ contains the set of
distinguished points of $S(G)$, but instead we suppose given a subset
$\Sigma\subset X_0(G)\sqcup X_{\infty}(G)$ which contains both $x(P)$ and the
distinguished points of $G$ and an
involution $\iota$ on the complement $\Sigma - x(P)$. We define inductively the
{\it order} of a connected component of $G$ as follows: a connected component
is of order zero if it contains a point of $x(P)\cap X_{\infty}(G)$; a
connected component has order $\le k+1$ if it contains a distinguished point
$p$ such that $\iota (p)$ lies on a component of order $\le k$.

We say that $(G,x,\iota )$ is a {\it stable $P$-pointed ribbon graph} if
\roster
\item every component has an order and
\item for every $p\in X_{\infty}(G)$ on a component of order $k>0$, $\iota (p)$
is on a component of order $k-1$.
\endroster
(So in the situation (2) we must have $\iota (p)\in X_0(G)$.)

\medskip\label
A stable $P$-pointed ribbon graph $(G,x,\iota )$
determines a stable $P$-pointed pseudosurface $(S(G,\iota ),x)$: it
is obtained from the surface $S(G)$ by identifying the points (of $\Sigma
-x(P)$) according to the involution $\iota $. If this surface is connected,
then it has a {\it genus $g$} characterized by the condition that $2-2g$ is the
Euler characteristic of the smooth part of $S(G,\iota )$.

We have seen that a conformal structure $l$ on $G$ determines a
conformal structure on $S(G)$ so that we have a compact Riemann surface $C(G
,l)$. This in turn, determines a weakly normal complex-analytic structure on
$S(G,\iota )$. With that structure, $(S(G,\iota ),x)$ becomes a stable
$P$-pointed
curve $(C(G,\iota ,l),x)$. This curve has additional structure: to every point
$p\in x(P)\cup S(G,\iota )_{\sing}$ is assigned a nonnegative
number $\lambda (p)$, namely half the length of the corresponding boundary
cycle (with respect to the componentwise unital metric defining the conformal
structure) in case the point comes from $X_{\infty}(G)$ and zero else. Notice
that $\lambda (p)=0$ if $x(p)$ lies on a single irreducible component of
$S(G ,\iota )$ or if $p\in P$ and $x(p)\in X_0(G)$, and that the sum of the
values of $\lambda$ on each irreducible component is $1$.

This  suggests to extend the notion of a {\it circumference function} to the
case of a stable connected $P$-pointed pseudosurface $(S',x)$ as as a function
$\lambda :x(P)\cup S' _{sing}\to\R _{>0}$  which possesses these properties.
So the space of circumference functions on $(S',x)$ is a product of simplices
(with a factor for each irreducible component).

\medskip\label
Just as for smooth $P$-pointed curves, the datum of a cicumference
function  $\lambda$ on a stable $P$-pointed curve $(C,x)$  permits us to go in
the opposite direction: apply Strebel's theorem
\refer{7.6} componentwise to the normalisation $(\hat C,\lambda )$. This
determines a Jenkins-Strebel differential $q$ on $\hat C$ with the properties
mentioned there. In particular, we have a critical graph $(G ,l)$ in $\hat C$
which contains the zeroes of $\lambda$. Moreover, each $p\in\supp (\lambda)$
determines (and is determined by) a boundary cycle of $G$ and the length of
that boundary
cycle is $\lambda (p)$. The associated Riemann surface $C(G ,l)$ is naturally
isomorphic to $\hat C$.

\medskip\label
Let now $(G,x)$ be a $P$-pointed ribbon graph. We describe how a proper subset
of $X_1(G)$ (or rather, strictly decreasing sequences of such) define stable
$P$-pointed ribbon graphs. First two definitions.

Let $Z$ be a semistable set of edges of $G$. Recall that then every component
of $G_Z$ that is a homotopy circle without a vertex in $x(P)$ is necessarily a
topological circle (and is not a boundary cycle of $G$). If this does not
happen, i.e., if every component of $G_Z$ that is a topological circle
contains a vertex in the image of $x$, then we say that $Z$ is {\it stable}. It
is clear that every subset $Z\subset X_1(G)$ has a maximal semistable subset
$Z^{\st}$; it is a union of components of $Z^{\sst}$.

Forgetting the bivalent vertices of $G_{Z^{\st}}$ that are in $x(P)$ yields
a ribbon graph with the same underlying topological space as $G_{Z^{\st}}$; we
denote this ribbon graph by $\bar G_{Z^{\st}}$ and its set of edges by $\bar
Z^{\st}$. It is clear that the set of distinguished points of $S(G_{Z^{\st}})$
coincides with $X_{\infty}(G_{\bar Z^{\st}})$.

A metric on $G_{Z^{\st}}$ determines one on $\bar G_{Z^{\st}}$.

\medskip\label Let $Z$ be a proper subset of $X_1(G)$ and put
$G(Z):=G/G_Z\sqcup G_{Z^{\st}}$. It is clear that the pointing $x$ determines
an injection $\tilde x$ of $P$ in the set of $0$-simplices of $G(Z)$.
The proof of the following lemma is easy and left to the reader

\proclaim{\label Lemma} The set of distinguished points of $G(Z)$ that are not
in the image of $\tilde x$ comes with a natural involution $\iota$ so that
$G(Z)$, $\tilde x$ and $\iota$ define a stable $P$-pointed ribbon graph. The
associated $P$-pointed stable pseudosurface $S(G;Z)$ is obtained from
$S(G/G_Z)$ and $G_{Z^{\st}}$ by identifying
each exceptional vertex of $S(G/G_Z)$ with the corresponding exceptional
element of $X_{\infty}(G_{Z^{\st}})$ and then contracting every irreducible
component that corresponds to a component of $G_{Z^{\sst}}-G_{Z^{\st}}$. A
conformal structure on $\tilde G$ determines one on $S(\tilde G,\iota)$ and
turns the latter into a stable $P$-pointed curve.
\endproclaim

\label We may of course repeat this construction for a set of edges of $G_{\bar
Z^{\st}}$. In order to be able to state this we introduce the following
notions.

A {\it permissible sequence} for $(G ,x)$ is a  sequence $Z _{\bullet}=(X_1(G )
\! =\! Z_0,Z_1, Z_2,\dots ,Z_k)$ such that $Z _{\kappa}\subset \bar Z_{\kappa
-1}^{\st}$ and $G_{Z _{\kappa}}$ does not contain a connected component of
$\bar
G_{Z_{\kappa -1}^{\st}}$.

A {\it stable metric} relative such a
sequence is given by a conformal structure on every difference
$\bar G _{Z_{\kappa }^{\st}}-G _{Z_{\kappa +1}}$. So this may be given  by a
sequence of functions $l_{\kappa}:Z_{\kappa}^{\st}\to\R {\ge 0} $ such that
$l_{\kappa}$ has zero set $Z_{\kappa +1}$ ($\kappa =0,1,\dots $). (So
$l_{\bullet}$ determines $Z _{\bullet}$.)

The previous discussion generalizes in a straightforward way to:

\proclaim{\label Proposition}
Let $Z _{\bullet}$ be a permissible sequence for $(G ,x)$. Then the disjoint
union of the ribbon graphs $G_{\bar Z^{\st}_{\kappa}}/G_{Z_{\kappa +1}}$
($\kappa =0,1,\dots )$  is in a natural way a stable $P$-pointed ribbon graph
$(G(Z_{\bullet}),\tilde x,\iota )$. A stable metric $l_{\bullet}$ relative
$Z _{\bullet}$ defines a conformal structure on  $S(G,Z_{\bullet})$ and turns
it into a stable $P$-pointed curve $C(G, l_{\bullet})$.
\endproclaim

\head
\section Stable limits
\endhead
In this section we fix a connected $P$-pointed ribbon graph $(G,x)$. We explain
how the stable pseudosurface associated to a permissible sequence for $G$
arises as a limit of Riemann surfaces $C(G,l(t))$.

\medskip\label
We shall use a blowing up construction in the PL-category.
The basic construction starts out from a collection $\beta$ of oriented edges
of
$G$ that defines an oriented circular subgraph $G _{\beta}$ of $G$. Let
$U_{\beta}$ be the union
of the relatively open simplices that have a  point of $G _{\beta}$ in their
closure;  this is a regular neighborhood of $G _{\beta}$ PL-homeomorphic to
an open cylinder. Notice that
$U_{\beta}-G _{\beta}$ has two connected components, one which contains
the interiors of the tiles associated to the elements of $\beta$; we denote
that component $U_{\beta}^+$ and the other by $U_{\beta}^-$. By means of the
barycentric coordinates of the simplices in $U^+_{\beta}$ we have defined a
piecewise-linear function $U^+_{\beta}\to [0,1)$ which measures the distance to
$G _{\beta}$. Let  $\phi _{\beta}: U_{\beta}\to [0,1)$ be its extension by zero
on $U_{\beta}$; this is a continuous PL-function. Let $(U_{\beta}\times
\R _{\ge 0})^{\,\widetilde{}}$ be the closure of the graph of the function
$$
(u,t)\in (U_{\beta}-G _{\beta})\times\R _{>0}\mapsto [-\log (1-\phi
_{\beta}(u):t]\in P^1(\R )
$$
in $U_{\beta}\times\R _{\ge 0}\times P^1(\R )$.
The projection  $(U_{\beta}\times \R _{\ge 0})^{\,\widetilde{}}\to
U_{\beta}\times\R _{\ge 0}$ is clearly a PL-homeomorphism over the complement
of $G
_{\beta}\times\{ 0\}$ whereas the pre-image of $G _{\beta}\times\{ 0\}$ is $G
_{\beta}\times\{ 0\}\times [0,\infty ]$. The strict transform of
$U_{\beta}^+\times\{ 0\}$ resp.\  $U_{\beta}^-\times\{ 0\}$ meets $G
_{\beta}\times\{ 0\}\times [0,\infty ]$ in  $G _{\beta}\times\{ 0\}\times
\{\infty \}$ resp.\ $G _{\beta}\times\{ 0\}\times \{ 0\}$. So the total
transform of $U_{\beta}\times\{ 0\}$ is a kind of thickening of $U_{\beta}$
(see
the figure below).

\midspace{50mm}\caption{\fig Blowing up of an oriented cycle}

In particular, this total transform is PL-homeomorphic
to $U_{\beta}$; indeed, the projection $(U_{\beta}\times
\R _{\ge 0})^{\,\widetilde{}}\to \R _{\ge 0}$ is trivial.

We glue $(U_{\beta}\times\R _{\ge 0})^{\,\widetilde{}}\to \R _{\ge 0}$ to
$(S(G )\times \R _{\ge 0})-(G _{\beta}\times\{ 0\})$ via their common open
subset
$U_{\beta}\times\R _{\ge 0} -G _{\beta}\times\{ 0\}$ and obtain a modification
$(S(G )\times\R _{\ge 0})_{\beta}^{\,\widetilde{}}\to S(G )\times\R _{\ge 0}$.

For $e\in\beta$, the tile $K_e\times\{ 0\}$ lifts PL-homeomorphically to
$(S(G )\times\R _{\ge 0})_{\beta}^{\,\widetilde{}}$. We apply an edge collapse
to all these lifted copies and denote the result $(S(G )\times\R _{\ge
0})_{\beta}^{\,\widehat{}}$.
The pre-image of $S(G)\times\{ 0\}$ is denoted by $S(G;\beta )$. It is a
pseudosurface that is PL-homeomorphic to the the space obtained from $S(G)$ by
contracting $G_{\beta}$. It comes with an injection of $P$ in its regular
part.

\medskip\label
We now fix a proper subset $Z$ of $X_1(G )$ and show how $S(G ;Z)$ is obtained
as a one-parameter degeneration of $S(G )$. First we assume that $Z$ is
stable.
We carry out the previous construction for each boundary cycle of $Z$. It is
easily seen that these can be performed independently so
that we have defined a modification
$$
(S(G )\times\R _{\ge 0})^{\,\widehat{}}_Z\to S(G )\times\R _{\ge 0}.
$$
The crucial remark is that this projection $(S(G )\times \R _{\ge
0})^{\,\widehat{}}_Z\to \R _{\ge 0}$ is trivial over $\R _{>0}$ with fiber $S(G
)$ whereas the fiber over $0$ is canonically isomorphic to $S(G;Z)$.

In case $Z$ is not stable, we first apply the preceding procedure to $Z^{\st}$
and next we collapse the strict transforms of the tiles indexed by the oriented
members of $Z-Z^{\st}$ (a total collapse or an edge collapse, depending on
whether the boundary cycle generated by the corresponding oriented edge is in
$G_Z$ or not). The order of these operations can be reversed; in particular, we
can first pass to $G':=G/G _{Z-Z^{\sst}}$ and the image $Z'$ in
$X_1(G')$ (so that $Z'$ is now semistable), then perform edge collapses on the
tiles indexed by the oriented members of $Z'-Z'{}^{\st}$ (these make up a union
of circular components of $G_{Z'}$) and finally apply the preceding
construction with $Z'{}^{\st}$.
Then the fiber over $0$ can be identified with $S(G;Z)$ as before.

\smallskip
We already observed that conformal structures $l_0$ on $G/G_Z$ and
$l_1$ on $G_Z$ determine a conformal structure on $S(G;Z)$, turning it into a
stable $P$-pointed curve $C(G,(l_0,l_1))$ whose normalization is the disjoint
union of the Riemann surfaces $C(G/G_Z,l_0)$ and $C(G_{Z^{\st}},l_1)$.
We may obtain such conformal structures  by means of a degeneration of a family
of metrics on $S(G)$. To be concrete, let $l$ be a metric on $G$ and let for
$t>0$, $l(t)$ be the metric on $G$ which takes on an edge $s$ the value $tl(s)$
if $s\in Z$ and remains $l(s)$ if not. We give the fiber of $(S(G )\times \R
_{\ge 0})^{\,\widehat{}}_Z\to \R _{\ge 0}$ over $t\in \R _{>0}$ (which is just
$S(G)$) the corresponding metric structure (denoted $m_t$).  The regular part
of
the fiber over $0$ is given the metric structure $m_0$ defined by the
restrictions $l_0$ resp.\ $l_1$ of $l$ to $X_1(G )-Z$ resp.\ $Z$. This is in
general not a continuous family of metrics, but for the underlying conformal
structures the situation is different. To see this, let $\phi _Z:S(G)\to \R
_{\ge 0}$ be the piecewise-linear function that takes the value $0$ on every
vertex in $G_Z$, $U_Z\subset S(G)$ the set where $\phi _Z<1$ and put
$f_Z:=-\log
(1-\phi _Z): U_Z\to\R _{\ge 0}$. It is clear from our definition of $m_t$ that
the set $f_Z<a$ with metric $m_t$ is conformally equivalent to subset
$f_Z<t^{-1}a$ with metric $m_1$. In fact, we have

\proclaim{\label Lemma}
Suppose that the pointing $x$ of $G$ has been extended
to a marking by $(S,P)$. Then the map $\R _{\ge 0}\to \overline{\T }$, which
assigns to
$t>0$ resp.\ $t=0$ the isomorphism class of $C(G,l(t))$ resp.\ $C(G,
(l_0,l_1))$ is continuous.
\endproclaim
\demo{Proof} There is no loss of generality in assuming that $G_Z$ has no
negligible components.

The continuity on $\R _{>0}$ is clear. To prove continuity at $0$ we wish to
invoke \refer{2.1}. This requires that we trivialize our family locally. At the
points of $S(G;Z)$ outside the exceptional set this is no problem and it is
clear that relative a suitable
trivialisation the complex structures converges uniformly on compact subsets.
At
the points of $S(G;Z)$ outside the strict transform we trivialize as follows.
Choose a piecewise-linear retraction $r_Z:U_Z\to G_Z$ so that $(r_Z,f_Z)$
defines a PL-homeomorphism $h$ of $U_Z-G_Z$ onto $\tilde G_Z\times\R_{>0}$,
where $\tilde G_Z$ is the disjoint union of the boundary cycles of $G_Z$. Let
$k$ denote its inverse and for $t>0$, let  $k_t(p,s)=k(p,st)$. Then
$$
(p,s,t)\in \tilde G_Z\times\R _{> 0}\times\R _{>0}\mapsto (k_t(p,s),t)
$$
extends to a PL-homeomorphism of $\tilde G_Z\times \R _{> 0}\times\R _{\ge 0}$
onto an open subset of $(S(G)\times \R _{\ge 0})^{\,\widehat{}}_Z$ so that for
$t=0$ we get a PL-homeomorphism $k_0$ of $\tilde G_Z\times \R _{> 0}$
onto the complement of the union of the strict transform of $S(G)$ and $G_Z$ in
$S(G;Z)$. We must show that the conformal
structure $J_t$, $t\ge 0$ on $\tilde G_Z\times (0,1)$ defined by pull-back of
the given conformal structure on $C(G,l(t))$ under $k_t$ depends continuously
on $t$. This is proved using explicit coordinates. We leave that to the
reader.
\enddemo

The preceding  can be iterated in an obvious way and yields:

\proclaim{\label Proposition}
If $Z_{\bullet}$ is a permissible sequence, then there is defined an iterated
modification:
$$
(S(G)\times \R _{\ge 0})^{\,\widehat{}}_{Z_{\bullet}}\to \R _{\ge 0}.
$$
This fibration is canonically trivialized (relative $x$) over $\R _{>0}$ with
fiber $S(G)$, whereas the fiber over $0$ is canonically homeomorphic to
$S(G;Z_{\bullet})$.

Suppose that the pointing $x$ of $G$ has been extended to a marking by
$(S,P)$.
Given a metric $l$ on $G$, let $l(t)$ be the metric on $G$ that
on $G_{Z_{\kappa}}-G_{Z_{\kappa +1}}$ is equal to $t^{\kappa}l$ ($t>0$) and let
$l_{\bullet}$ be the stable metric relative  $Z _{\bullet}$ that is defined
by the restrictions of $l$. Then the map $\R_{\ge 0}\to \overline{\T }$
which assigns to $t\in\R_{>0}$ resp.\ $t=0$ the isomorphism type of $C(G,l(t))$
resp.\ $S(G)\cong C(G,l_{\bullet})$ is continuous.
\endproclaim

\head
\section Deligne--Mumford modification of the arc complex
\endhead

Let $(G,x)$ be a connected $P$-pointed ribbon graph. Recall that we have
defined
the family $\pi :\CC (G )\to a(G)$ which over the interior $\conf (G )$ of
$a(G)$ is trivialized with fiber $S(G)$.  We are going to modify this family
over
the locus where this family is not locally trivial. This will also modify the
base and the result will be a family parametrizing stable pointed
pseudosurfaces  with stable conformal structures.

\medskip
Let $\Cal{Z}(G )$ denote the collection of stable subsets $Z\subset X_1(G )$
with $G_Z$ connected. For $l\in \conf (G )$ and $Z\in \Cal{Z}(G )$, we let $\pi
_Z(l)$ denote the unital metric on $G _{\bar Z}$ which is proportional to $l|G
_Z$. Let $\hat a(G )$ be the closure of the graph of the map $l\in \conf (G
)\mapsto (\pi _Z(l)\in \conf (G _Z))_Z$ in $a(G )\times \prod _{Z\in
\Cal{Z}(G
)} a(G _{\bar Z})$.

\proclaim{\label Proposition}
There is a natural bijection between the points of $\hat a(G )$ and the set of
stable conformal structures on $G$.
\endproclaim

\demo{Proof}
Let $(l^{(n)})_{n=1}^{\infty}$ be a sequence in $\conf (G )$. By passing to a
subsequence, we may assume
that for every $Z\in\ZZ (G)$, the sequence $(\pi _Z(l^{(n)}))_n$ converges
(to $l_Z$, say).  Write $l_0$ for $l_{X_1(G)}$, let $Z(l_0)$ be the zero set of
$l_0$ and put $Z_1:=Z(l_0)^{\st}$. Notice that
$\Cal{Z}(G_{Z_1})$ is just a subset of $\Cal{Z}(G_Z)$. So for each $Z\in
\Cal{Z}(G_{Z_1})$  we have a function $l_{Z}:Z\to [0,1]$ whose sum is $1$.
Applying this to the connected components of $G_{Z_1}$ yields a function
$l_1:Z_1\to [0,1]$ that on each connected component of $G_{Z_1}$ sums up to
$1$. We proceed with induction: if $l_{\kappa}: Z_{\kappa}\to [0,1]$ has been
constructed, then let $Z(l_{\kappa})$ be the zero set of $l_{\kappa}$. We put
$Z_{\kappa +1}:=Z(l_{\kappa})^{\st}$ and define $l_{\kappa +1}: Z_{\kappa
+1}\to [0,1]$ by letting it on each connected component $G_{Z}$ of
$G_{Z_{\kappa +1}}$ be equal to $l_{Z}$. Then $Z_{\bullet}$ is a permissible
sequence for $(G ,x)$ by construction. It comes naturally with a unital stable
metric $l_{\bullet}$ relative this sequence. This stable metric determines
every $l_Z$: for $Z\in\Cal{Z}(G)$, let $\kappa$ be such that $G_Z\subset G_{Z
_{\kappa}}$ and $G_Z\not\subset G_{Z _{\kappa +1}}$. Then  $G_Z$ is contained
in
a connected component $G_{Z'}$ of $G_{Z_{\kappa}}$. Since $Z\not\subset Z
_{\kappa +1}$, $l_{Z'}|Z$ (and hence $l_{\kappa }|Z$) is not identically zero.
It then follows that $l_Z$ is the unital metric proportional to $l_{\kappa
}|Z$.
On the other hand, \refer{10.4} shows that every  stable metric thus arises.
\enddemo

\medskip\label
If $Z$ is a negligible set of edges of $G$, then $\hat a(G /G _Z)$ can be
identified with the subset of $\hat a(G )$ parametrizing stable metrics
$l_{\bullet}$ of which each term vanishes on $Z$. Hence if we endow the ribbon
graphs with markings, then the closed cells $\hat a(G ,f)$ can be glued
together to yield a modification
$$
\hat A\to A.
$$
It is clear that $\hat A$ comes with a decomposition into cells. Such a cell
admits a description in terms of arc complexes as follows:
it is of the form $\sigma _0\times\sigma _1\times\dots$, where each $\sigma
_{\kappa}$ is a cell (a product of simplices)  of the arc complex associated to
a (not necessarily connected) pointed surface $(S_{\kappa},P_{\kappa})$.
These pointed surfaces (and hence these
cells) are defined inductively: $(S_0,P_0):=(S,P)$ and $\sigma _0$ is an
arbitrary simplex of $A$.
For $\kappa\ge 1$, let $\bar S'_{\kappa}$ be the pseudosurface  obtained from
$S_{\kappa -1}$ by contracting the arcs that make up $\sigma _{\kappa -1}$,
$S'_{\kappa}$ its normalisation, and let $P'_{\kappa}\subset S'_{\kappa}$ the
pre-image of the image of $P_{\kappa -1}$. Let $(S_{\kappa},P_{\kappa})$ be
obtained from $(S'_{\kappa},P'_{\kappa})$ by discarding all components that are
one- or two-pointed spheres. The connected components of $S_{\kappa}$ label the
factors of $\sigma _{\kappa}$ so that each factor is made up of arcs in that
component. We require that these arcs connect only points of $P_{\kappa}$ that
map to singular points of $\bar S'_{\kappa -1}$.
Under the projection $\hat A\to A$ this cell maps to $\sigma _0$.

It is possible to give a complete description of the incidence relations
between these cells, but we omit this.

\medskip\label
We shall define a family of surfaces $\hat\CC (G )$ over
$\hat a(G )$.  Let $Z_{\bullet}$
be a permissible sequence for $G$ of {\it connected stable} subsets, which we
here  regard as a strictly decreasing sequence of connected stable subsets of
$X_1(G)$, and consider the map
$$
I_{Z_{\bullet}}: S(G)\times\conf (G)\to \prod _{\kappa \ge 1}
(S(G)\times \R _{>0}),\quad
(u,l)\mapsto (u,l(Z_{\kappa})/l(Z_{\kappa -1}))_{\kappa}.
$$
The closure of its graph in
$S(G)\times a(G)\times \prod _{\kappa \ge 1}(S(G)\times \R _{\ge
0})^{\,\widehat{}}_{Z_{\kappa}}$ is denoted  by
$(S(G)\times a(G))^{\,\widehat{}}_{Z_{\bullet}}$.

Similarly, we denote the closure of the graph of
$$
\conf (G)\to\prod _{\kappa \ge 1}\R _{>0},\quad l\mapsto
(l(Z_{\kappa})/l(Z_{\kappa -1}))_{\kappa}
$$
in $a(G)\times\prod _{\kappa \ge 1}\R _{\ge 0}$ by
$a(G)^{\,\widehat{}}_{Z_{\bullet}}$. Since the functions
$l(Z_{\kappa})/l(Z_{\kappa -1})$ extend continuously to $\hat a(G)$, this is a
quotient of $\hat a(G)$.
We have a projection
$$
(S(G)\times a(G))^{\,\widehat{}}_{Z_{\bullet}}\to
a(G)^{\,\widehat{}}_{Z_{\bullet}}.
$$
Any fiber over a point of $a(G)^{\,\widehat{}}_{Z_{\bullet}}$ that has all its
coordinates in $\prod _{\kappa \ge 1}\R _{\ge 0}$ equal to zero is isomorphic
to $S(G;Z_{\bullet})$.

\smallskip
We do this for all such sequences simultaneously. To be precise, let $\ZZZ (G)$
be the collection of strictly decreasing sequences of connected stable subsets
of $X_1(G)$, and consider the map
$$
I=(I_{Z_{\bullet}}):S(G)\times\conf (G)\to \prod _{Z_{\bullet}} \prod _{\kappa
\ge 1} (S(G)\times \R _{> 0}).
$$
The closure of its graph in
$$
\hat a(G)\times \prod _{Z_{\bullet}} \prod _{\kappa \ge 1} (S(G)\times \R _{\ge
0})^{\,\widehat{}}_{Z_{\kappa}}
$$
is denoted $\hat\CC (G)$ and the projection of $\hat\CC (G)$ onto $\hat a(G)$
by $\hat \pi _G$. The preceding discussion shows:

\proclaim{\label Proposition}
If $l_{\bullet}$ is a stable metric with associated
permissible sequence $Z_{\bullet}$, then the fibre $\hat \pi _{G }
^{-1}(l_{\bullet})$ is naturally homeomorphic to $S(G;Z_{\bullet})$.
\endproclaim

We endow the fiber $\hat \pi _{G }^{-1}(l_{\bullet})$ with the conformal
structure prescribed by the stable metric $l_{\bullet}$ so that $\hat\pi _G$
defines a family of stable $P$-pointed stable curves.

For marked ribbon graphs this construction is
compatible in the sense that if $Z\subset X_1(G )$ is negligible, then  $\hat
\pi _{G /G _Z}:\hat\CC (G /G _Z)\to \hat a(G /G _Z)$ can be identified with
the restriction of $\hat \pi _{G }$ over $\hat a(G /G _Z)$. We may therefore
glue these maps to each other to get a family  $\hat \pi :\hat\CC \to \hat
A$ of stable $P$-pointed curves.  Each fiber of $\hat \pi $ maps to a fiber of
$\pi $, so that we have a commutative diagram
$$ \CD
\hat\CC @>>> \CC \\
@V\hat \pi VV        @V\pi VV\\
\hat A @>>> A
\endCD
$$
of spaces with $\G$-action. We have also have a classifying map that
extends $\Phi$:
$$
\hat\Phi :\hat A\to\overline{\T } .
$$
It is clearly $\G$-equivariant. Our second main result reads as follows:

\proclaim{\label Theorem}
The map $\hat\Phi:\hat A\to \overline{\T }$ is a $\G $-equivariant continuous
surjection. The pre-image of the class of a marked stable  $P$-pointed curve
$(C,[f])$ under $\hat\Phi$ can be
identified with the space of circumference functions \refer{9.2} of $(C,x)$. In
particular,
$\hat\Phi$ drops to a continuous surjection of $\G \backslash\hat A $ onto the
Deligne--Mumford compactification $\overline{\Mod}{}^P_g$.
\endproclaim

\demo{Proof} Let $(C,[f])$ be as in the theorem. The construction described in
\refer{9.3} produces for every circumference function of $(C,x)$ a marked
ribbon graph $(G ,f)$ plus a stable metric $l_{\bullet}$ on $G$ which
reconstructs $(C,[f])$ for us. This defines an element of $\hat a(G ,f)$ and
one verifies that its image in $\hat A$ is unique.

It remains to show that $\hat\Phi$ is continuous. It is enough to prove that
its
restriction to every closed cell $\hat a(G ,f)$ is. Since $\hat a(G ,f)$ is
second countable and $\overline{\T }$ is Hausdorff, we only need to verify
that the image of a converging sequence $(l_{\bullet}^{(n)})_n$ in $\hat a(G
,f)$ under $\hat\Phi$ has a limit point. Then after
passing to a  subsequence we may assume that $(l_{\bullet}^{(n)})_n$ is in the
relative interior of a single cell, say of $\hat a(G ,f)$. The desired property
then follows from \refer{2.1} as in the proof of \refer{10.3}.
\enddemo

We can now finish the proof of our first main result, too.

\medskip
\demo{Proof of \refer{8.6}}
The map $\hat\Phi$ and the projection
$\hat A\to A\to\Delta _P$ together define a map from $\hat A$ to
$\overline{\T }\times\Delta _P$. If we compose the latter with the quotient
map $\overline{\T }\times\Delta _P\to |K_{\bullet}\T |$ we get a map
$\hat\Psi :\hat A\to |K_{\bullet}\T |$. The theorem above implies that the
fibers of $\hat\Psi$ and the fibers of $\hat A\to A$ coincide. The
induced bijection $A\to |K_{\bullet}\T |$ is just $\Psi$. Since
$A$ has the quotient topology for the projection $\hat A\to A$,
it follows that $\Psi$ is continuous.  \enddemo


\Refs

\ref\no 1
\by B.H.\ Bowditch \&\ D.B.A.\ Epstein
\paper Natural triangulations associated to a surface
\jour Topology
\vol 27
\pages 91-117
\yr 1988
\endref

\ref\no 2
\by P.\ Deligne \&\ D.\ Mumford
\paper The irreducibility of the space of curves of given genus
\jour Inst\. Hautes \'Etudes Sci\. Publ\. Math\.
\vol 36
\yr 1969
\pages 75--109
\endref


\ref\no 3
\by J.L.\ Harer
\paper The virtual cohomological dimension of the mapping class group
\jour Invent. Math.
\vol 84
\yr 1986
\pages 157--176
\endref

\ref\no 4
\by J.L.\ Harer
\paper The  cohomology of the moduli space of the space of curves
\pages 138--221
\inbook Theory of Moduli
\ed E.\ Sernesi
\bookinfo Lecture Notes in Math.
\vol 1337
\publ Springer
\publaddr Berlin and New York
\yr 1988
\endref

\ref\no 5
\by W.J.\ Harvey
\paper Boundary structure of the modular group
\inbook Riemann surfaces and related topics
\bookinfo Annals of Math. Studies
\eds I.\ Kra and B.\ Maskit
\publ Princeton UP
\yr 1981
\pages 245--251
\endref

\ref\no 6
\by M.\ Kontsevich
\paper Intersection theory on the moduli space of curves and the
matrix Airy function
\jour Comm. Math. Phys.
\vol 147
\yr 1992
\pages 1--23
\endref


\ref\no 7
\by J.\ Milgram \&\ R.C.\ Penner
\paper Riemann's moduli space and the symmetric group
\inbook Mapping class groups and moduli spaces of Riemann surfaces
\eds C.--F.\ B\"odigheimer \&\ R.\ Hain
\bookinfo Contemp.\ Math.
\vol 150
\pages 247--290
\yr 1993
\endref


\ref\no 8
\by R.C.\ Penner
\paper Perturbative series and the moduli space of punctured surfaces
\jour J.\ Diff.\ Geom.
\vol 27
\yr 1988
\pages 35--53
\endref

\ref\no 9
\by K.\ Strebel
\book Quadratic differentials
\bookinfo Ergebnisse der Math.\ u.\ ihrer Grenzgebiete, 3.\ Folge
\vol 5
\publ Springer
\publaddr Berlin and New York
\yr 1984
\endref

\ref\no 10
\by S.\ Wolpert
\paper The hyperbolic metric and the geometry of the universal curve
\jour J.\ Diff.\ Geom.
\vol 31
\pages 417--472
\yr 1990
\endref

\endRefs
\enddocument